\documentclass[11pt]{article}%
\usepackage{amssymb}
\usepackage{amsfonts}
\usepackage{amsmath}
\usepackage{amssymb}
\usepackage{color}
\usepackage{graphicx}
\usepackage{caption}%
\usepackage[authoryear]{natbib}
\usepackage{hyperref}
\usepackage{float}
\usepackage{subcaption}
\usepackage{multirow}

\usepackage{tabularx}
\usepackage{algorithm}
\usepackage{algpseudocode}
\providecommand{\U}[1]{\protect\rule{.1in}{.1in}}
\newtheorem{theorem}{Result}

\begin{document}
	
	\title{\textbf{Robust inference for an interval-monitored step-stress experiment with competing risks for failure}}
	\author{Narayanaswamy Balakrishnan,  María Jaenada and Leandro Pardo}
	\date{ }
	\maketitle

\begin{abstract}
	Accelerated life-tests (ALTs) are used for inferring  lifetime characteristics of highly reliable products.
	In particular, step-stress ALTs increase the stress level at which units under test are subject at certain pre-fixed times, thus accelerating the product's wear and inducing its failure. In some cases, due to cost or product nature constraints, continuous monitoring of devices is infeasible, and so the units are inspected for failures at particular inspection time points. In a such setup, the ALT response is interval-censored.
	Furthermore, when a test unit fails, there are often more than one fatal cause for the failure, known as competing risks. In this paper, we assume that all competing risks are independent and follow  exponential distributions with scale parameters depending on the stress level.
	Under this setup, we present a family of robust estimators based on density power divergence, including the classical maximum likelihood estimator (MLE) as a particular case.
	We derive asymptotic and robustness properties of the Minimum Density Power Divergence Estimator (MDPDE), showing its consistency for large samples.
	Based on these MDPDEs, estimates of the lifetime characteristics of the product as well as estimates of cause-specific lifetime characteristics  are then developed. Direct asymptotic, transformed and bootstrap confidence intervals for the mean lifetime to failure, reliability at a mission time and distribution quantiles are proposed, and their performance is then compared through Monte Carlo simulations.
	Moreover, the performance of the MDPDE family has been examined through an extensive numerical study and the methods of inference discussed here are  finally illustrated with a real-data example concerning electronic devices.

\end{abstract}

\textbf{Keywords}:
	Accelerated Life-Tests, Competing Risks, Minimum Density Power Divergence Estimator, Reliability Tests, Robustness.

\section{Introduction} \label{sec:Introduction}

Due to advancements in manufacturing design and technology, modern products have become increasingly reliable.
However, this progress presents challenges when it comes to life-testing  of these products. Testing the products under normal operating conditions would result in  few failures if at all, which would be insufficient for accurate inference. Moreover, life-tests under normal operating conditions would require long experimental time and consequently a high experimental cost. 
Therefore, conventional life-testing experiments under normal operating conditions may be unsuitable for highly reliable products.

Alternatively, Accelerated Life Testing (ALT) experiments would be more effective, producing more failures in lesser experimental times.
ALTs induce failures by increasing one or more environmental stress factors affecting the lifetime of the device, such as temperature, voltage, or humidity. After analyzing the data for the increased stress levels, results can be extrapolated to normal operating conditions.
There exist three main types of ALTs depending on how the increased stress is applied to the units under test. 
In constant-stress ALTs, test units are exposed to a fixed and elevated stress level throughout the duration of the test.
Therefore, the units are only subjected to a single stress level, but that level may change for different groups of units.
In contrast, progressive ALTs allow the stress to increase linearly and continuously on any surviving test unit. 
On the other hand, the step-stress tests allow an increase of the stress level gradually at some pre-planned time points, referred to here as times of stress change. Between two subsequent times of stress change, the stress level is maintained at a constant level.
That is, all units under test experience the same stresses, which are increased to induce more failures.
If the first stress level considered is the normal operating level, then the loading scheme applied is referred to as partial step-stress ALT.
Step-stress models have been widely studied in the literature as they offer a suitable trade-off between experimental flexibility and accuracy in inference.

Although the step-stress ALT models usually refer to continuous monitoring, in some experiments, the exact failure times may not be available and only the number of failures observed in specific time intervals may be recorded. In such a setup, the ALT response is interval-censored and the information available is the count of failures in specific time intervals. 
As the available information may be different for different types of censoring, specific statistical techniques must be developed for each censoring setup. 
Here, we will assume that the data are interval-censored.
This type of interval sampling is also described in the literature under the term ``grouped samples'', ``interval-monitoring data'' or ``non-destructive one-shot devices data''.

Competing risks arise in experiments in which units are subject to several potential failure causes and the occurrence of one event might impede the occurrence of other events.
This experimental model is observed in various fields, including engineering experiments and clinical analysis.
In engineering experiments, when investigating the lifetime to failure of a mechanical component like a gearbox, competing risks can emerge. For instance, the gearbox may have two competing failure modes: bearing failure and gear tooth failure. The occurrence of one failure mode can prevent the other from happening. Analyzing the competing risks helps engineers understand the reliability characteristics of the component and optimize its design and operating conditions.
In clinical analysis, competing risks are also prevalent. Consider the example of a patient who undergoes a bone marrow transplant. After the transplant, the patient faces the possibility of experiencing a relapse or passing away while in remission. In this case, relapse and death are competing risks because once the patient dies, the event of relapse becomes impossible. Studying these competing risks in clinical analysis helps in understanding the outcomes and survival probabilities of patients, enabling better treatment and health care decisions.
%
In the mentioned examples, the different failure causes of the devices can be described by a competing risks model, where in each risk represents a different cause of failure. In this context, the word "competing" refers to the fact that these failure causes are mutually exclusive, meaning that the occurrence of one failure mode prevents the occurrence of the other.
Also, it is assumed that the competing risks are independent of each other, and so the damage caused by one risk does not affect the reliability of the product due to other risk.
There are various examples of failure data under competing risks in the literature. 
\cite{balakrishnan2008exact} developed exact inference for a simple step-stress model with competing risks for exponential lifetimes under Type-II censoring. 
Along similar lines, \cite{han2010inference} developed exact inference for a simple step-stress model for the case of the exponential lifetime distribution, but under time constraints.
\cite{han2014inference} developed inferential methods for a step-stress model with competing risks assuming generalized exponential distributions under Type-I censoring.
More recently, \cite{mondal2022robust} proposed divergence-based estimation for dependant competing risks under interval monitoring for reliability tests under normal operating conditions.
For extreme interval-censoring schemes, 
\cite{balakrishnan2015em} studied MLEs for one-shot device testing with competing risks under exponential distribution and \cite{balakrishnan2023power} proposed robust divergence-based inference for the same set-up.

Most of the above-mentioned inferential methods use the MLEs for estimating the model parameters because of its well-known properties such as asymptotic efficiency, consistency or sufficiency. However, the MLEs may get heavily affected by outlying observations. As an alternative, recent works on step-stress models have shown the advantage of using divergence-based estimators for robust estimation;
\cite{balakrishnan2022robust, balakrishnan2023step, balakrishnan2023Non} used
the  MDPDEs for step-stress model under interval-censoring and exponential, gamma and lognormal distribution, respectively. The MDPDEs demonstrated there an appealing balance between efficiency and robustness in the estimation process. Here, we develop inferential methods based on the density power divergence (DPD) for simple step-stress ALTs with competing risks under exponential lifetime distributions.

\section{The step-stress model under competing risks \label{sec:SSALT}}


Let us consider a simple step-stress ALT with two stress levels $x_1$ and $x_2,$ $R$ independent competing risks causing failure and $N$ units under test. We denote by $\tau_1$ the time of stress change from $x_1$ to $x_2$ and by $\tau_2$ the experiment termination time. 
Moreover, we consider $L$ pre-fixed inspection times at which failure counts are recorded, $ 0 = IT_{0} < IT_1 < \cdots < IT_L = \tau_2.$ We assume that these inspection times include the time of stress change $\tau_1.$

Since step-stress tests increase the stress level at some time, a model relating the lifetime distribution of experimental units at one
stress level to the distributions at preceding stress levels is necessary for inference. 
Three main statistical models have been proposed in the literature for  assessing the effect of increased stress levels. The
tampered random variable model (TRV) proposed by \cite{degroot1979bayesian} scales down the remaining lifetime at successive stress levels. 
On the other hand, the tampered failure rate model (TFR) proposed by \cite{bhattacharyya1989tampered} assumes that the effect of change of stress is to multiply the initial failure rate function by a factor after the stress change time. 
Finally, the cumulative exposure model (CEM), first studied by \cite{Sedyakin1966} and discussed further by  \cite{Bagdonavicius1978} and \cite{nelson1990accelerated}, assumes that the residual life of the experimental units depends only on the cumulative exposure the units have experienced, with no memory of how this exposure was accumulated. 
The CEM is one of the most commonly used models in step-stress inference and so we will adopt the CEM approach here. For a detailed review of maximum likelihood inferential methods for exponential step-stress models, one may refer to \cite{balakrishnan2009synthesis}.

Specifically, let us denote $F^{(1)}$ and $F^{(2)}$ for the cumulative distribution functions (CDFs) of the lifetimes of devices under the constant stresses $x_1$ and $x_{2},$ respectively. The CEM states that the damage incurred by the product until the time of stress change $\tau_1$ is embodied as a shift in the distribution $F^{(2)},$ such that the overall CDF of the lifetime is always continuous. That is, the shifting time $h$ must satisfy the equation
$$F^{(1)}(\tau_1) = F^{(2)}(\tau_1 + h).$$
Note that if the lifetime CDFs $F^{(i)},$ $i=1,2,$ belong to the same scale family of distributions, such as the exponential family,  with scale parameter $\theta_i$, the shifting time can be explicitly obtained as
 $$ h = \frac{\theta_2}{\theta_1}\tau_1 - \tau_1.$$


Assuming exponential distributions for the lifetimes of devices for each competing risk $T_j$, where $j=1,..,R$, at a constant stress level $x_i$, for $i=1,2$, with scale parameter $\theta_{ij} > 0$, the CDF of the lifetime attributed to cause $j$ can be expressed as 
$$F_j(t) = \begin{cases}
	1-\exp(-\frac{t +h^{(1)}_j}{\theta_{1j}}), & 0 < t < \tau_1,\\
	1- \exp(-\frac{t+h^{(2)}_j}{\theta_{2j}}), & \tau_1 \leq t < \infty,
\end{cases}
$$
with
\begin{equation} \label{eq:shiftingtime}
	h^{(1)}_j = 0 \hspace{0.3cm} \text{and} \hspace{0.3cm} h^{(2)}_j = \frac{\theta_{2j}}{\theta_{1j}}\tau_1 - \tau_1 
\end{equation}  for all $j=1,...,R.$ Here, we have included the shifting time at the first step for notational simplicity in the following expressions.
The corresponding probability density function (PDF) of $T_j$ is given by
$$f_j(t) = 
	\frac{1}{\theta_{ij}}\exp \left(-\frac{t+h^{(i)}_j}{\theta_{ij}} \right) \text{ if } \tau_{i-1} \leq t < \tau_{i},
$$
with $\tau_0 = 0$ and $\tau_2$ can be extended toward infinity.
Note that the shifting times are different for each competing risk, and depend on the distribution parameters of each lifetime $T_j.$

Since the risks are mutually exclusive, we will observe only the smaller time of failure $T=\operatorname{min}_{j=1,...,R}(T_j).$ Then, the CDF of the overall failure time of a test unit is readily obtained as
\begin{equation}
	F_T(t) = 1 - \prod_{j=1}^R(1-F_j(t)) 
		 =
		1- \exp\left(- \sum_{j=1}^R \frac{t+h_j^{(i)}}{\theta_{ij}}\right), \hspace{0.3cm} \tau_{i-1} \leq t < \tau_{i},
\end{equation}
with $h^{(i)}_j$ being as defined in Equation (\ref{eq:shiftingtime}), for any $i=1,2$ and $j=1,...,R.$
The corresponding PDF is consequently given by
\begin{equation}\label{eq:PDFmarginal}
		f_T(t) 
		=
			\left[\sum_{j=1}^R \frac{1}{\theta_{ij}}\right] \exp\left(- \sum_{j=1}^R \frac{t+h_j^{(i)}}{\theta_{ij}}\right), \hspace{0.3cm} \tau_{i-1} \leq t < \tau_{i}.
\end{equation}

Furthermore, let $C$ denote the indicator for the cause of failure. Then, the joint PDF of the failure time $T$ and the indicator failure mode $C$ can be expressed using the marginal distributions as
\begin{equation} \label{eq:jointPDF}
	f_{(T,C)}(t, j) = f_j(t)\prod_{\substack{j^\ast =1\\j^\ast \neq j}}^{R} (1-F_j^\ast(t))
	= 
			\frac{1}{\theta_{ij}} \exp\left(- \sum_{j=1}^R\frac{t+h^{(i)}_j}{\theta_{ij}}\right), \hspace{0.3cm} \tau_{i-1} \leq t < \tau_{i}.
\end{equation}
%
From the above joint density function, we can compute conditional probabilities, such as the relative risk imposed on a test unit before the time of stress change $\tau_1$ due to the risk factor $j,$ as
\begin{equation} \label{eq:pi1}
	\pi_{1j}(\boldsymbol{\theta}) = P\left( C=j | 0 < T < \tau_1 \right) = \frac{\theta_{1j}^{-1}}{\sum_{j=1}^R \theta_{1j}^{-1}}
\end{equation} 
 and similarly, the relative risk after $\tau_1$ due to the factor $j,$ as
\begin{equation} \label{eq:pi2}
	\pi_{2j}(\boldsymbol{\theta}) = P\left( C=j | \tau_1 < T < \infty \right)  = \frac{\theta_{2j}^{-1}}{\sum_{j=1}^R \theta_{2j}^{-1}} .
\end{equation} 
where $\boldsymbol{\theta} = \left(\theta_{11}, ...,\theta_{1R},\theta_{21} ...., \theta_{2R}\right).$
Both relative risks are the proportion of failure rates in the given time frame.

As we are dealing with interval-monitoring experiments, the resulting data would be failure counts within each inspected interval $(IT_{l-1}, IT_{l}],$ for $l=1,...,L.$ 
From the joint PDF in (\ref{eq:jointPDF}), we can compute the theoretical probability of failure within an interval due to each competing risk, as
\begin{equation} \label{eq:th}
	\begin{aligned}
	p_{lj}(\boldsymbol{\theta}) &= P\left(IT_{l-1} < t < IT_{l}, C=j \right) = \int_{IT_{l-1}}^{IT{l}} f_{(T,C)} (t,C=j) dt\\
	 &= 
	 \pi_{ij}(\boldsymbol{\theta}) \left(\exp\left(-\sum_{j=1}^{R} \frac{IT_{l-1}+h_j^{(i)}}{\theta_{ij}}\right) - \exp\left(-\sum_{j=1}^{R} \frac{IT_{l}+h^{(i)}_j}{\theta_{ij}} \right) \right), & \tau_{i-1} \leq IT_l < \tau_i.\\
	\end{aligned}
\end{equation}
Note that  $i$ corresponds to the stress level at which the units are subjected  during the period from $IT_{i-1}$ to $IT_{i}$ and $\pi_{ij}(\boldsymbol{\theta})$ is the relative risk imposed on a test unit before the time of stress change due to the risk factor $j,$ as defined in Equations (\ref{eq:pi1}) and (\ref{eq:pi2}), for $i=1$ and $i=2,$ respectively. Furthermore, we can assume that the stress level is constant within all intervals because the time of stress $\tau_1$ is included in the inspection times and the CDF is continuous in $\tau_1.$
On the other hand, the probability of survival is given by
$$p_0(\boldsymbol{\theta}) = 1 - F_T(\tau_2) = \exp\left(- \sum_{j=1}^R \frac{\tau_2+h_j^{(2)}}{\theta_{2j}} \right).$$

Because we have considered all exclusive failure events, all defined probabilities sum up to one,
$$\sum_{j=1}^R \sum_{l=1}^{L} p_{lj}(\boldsymbol{\theta}) + p_0(\boldsymbol{\theta}) = 1.$$

Given that the observed data consist of counts of events, the previous equation provides an understanding of considering a multinomial model with $N$ trials (devices) and $R \cdot L + 1$ events; the first $R \cdot L$ events correspond to failures attributed to each competing risk at each inspection interval, while the last event represents the survival beyond the experiment's termination.
The probability vector of the multinomial model is then given by $\boldsymbol{p}(\boldsymbol{\theta}) = (p_{11}(\boldsymbol{\theta}),..., p_{1R}(\boldsymbol{\theta}), p_{21}(\boldsymbol{\theta}),..., p_{2R}(\boldsymbol{\theta}), ....,  p_{L1}(\boldsymbol{\theta}),
 p_{LR}(\boldsymbol{\theta}), p_0(\boldsymbol{\theta})).$

Finally, a common assumption in ALT modeling is that the scale parameter of each marginal distribution is related to the stress level through a log-linear relationship. That is, for any competing risk $j=1,...,R,$ and stress level $x_i,$ $i=1,2$, we have $$\theta_{ij} = \exp(a_{0j} + a_{j2}x_i)$$
and so the model can be re-parametrized in terms of the parameter vector $\boldsymbol{a} = (\boldsymbol{a}_j)_{j=1,...,R}$ with $\boldsymbol{a}_j = (a_{0j}, a_{1j}).$
The log-linear relationship fits (upon suitable transformations) physical models relating stress factors and the lifetime of products, such as the Arrhenius law or the Inverse Power relationship, and for multiple step-stress ALTs helps in reducing the number of model parameters.

In terms of the regression parameter $\boldsymbol{a},$ we can rewrite the probability of failure at the $l$-th interval due to risk $j$ as

\begin{equation}
	\begin{aligned}
		p_{lj}(\boldsymbol{a})  = 
		\frac{\exp(-a_{0j}-a_{1j}x_i)}{\sum_{j=1}^R \exp(-a_{0j}-a_{1j}x_i)} 
		& \left(\exp\left(-\sum_{j=1}^{R} \frac{IT_{l-1}+h_j^{(i)}}{\exp(a_{0j}+a_{1j}x_i)}\right) \right.\\
		 & \left. - \exp\left(-\sum_{j=1}^{R} \frac{IT_{l}+h^{(i)}_j}{\exp(a_{0j}+a_{1j}x_i)} \right) \right), \hspace{0.3cm} \tau_{i-1} \leq IT_l < \tau_i,\\
	\end{aligned}	
\end{equation}
 and the probability of survival as
$$p_0(\boldsymbol{a}) = \exp\left(- \sum_{j=1}^R \frac{\tau_2+h_j^{(2)}}{\exp(a_{0j} + a_{1j}x_2)} \right),$$
with
\begin{equation} \label{eq:shiftingtime2}
	h^{(1)}_j = 0 \hspace{0.3cm} \text{and} \hspace{0.3cm} h^{(2)}_j = \tau_1\left(\frac{\exp(a_{0j}+a_{1j}x_2)}{\exp(a_{0j}+a_{1j}x_1)}- 1\right).
\end{equation} 
 
\section{The minimum density power divergence estimator}

The step-stress ALT model with competing risks described in Section \ref{sec:SSALT} relies on a $2R$-dimensional parameter vector $\boldsymbol{a}.$ Therefore, estimates of the model parameters  are required for analyzing the reliability of a device. 
The MLE is the most commonly used estimator for general statistical models. It possesses desirable properties including asymptotic efficiency and consistency. However, it lacks robustness and so a small amount of contamination in the observed data could heavily influence the estimation.
To overcome the robustness drawback, divergence-based estimators can be adopted. Specifically, DPD-based estimators have shown a gain in robustness with a small loss of efficiency in absence of contamination with respect to the MLE in step-stress models (\cite{balakrishnan2022robust, balakrishnan2022restricted, balakrishnan2023step}).
We develop here the MDPDE for step-stress ALT models with competing risks.

Let us first obtain the MLE for the step-stress model with competing risks and interval-monitoring.
We denote by $n_{lj}$ the number of failures recorded at the $l$-inspection time due to risk $j,$ $l=1,...,L$ and $j=1,...,R.$
Adopting the multinomial model discussed in Section \ref{sec:SSALT} with probability vector $\boldsymbol{p}(\boldsymbol{a})$, the (incomplete) log-likelihood of the step-stress model is given by
$$\mathcal{L}(\boldsymbol{a} | n_{11},...,n_{1R},...n_{LR}, n_0) \propto \sum_{j=1}^R\sum_{l=1}^{L} n_{lj}\log(p_{lj}(\boldsymbol{a})) + n_0\log(p_0(\boldsymbol{a})) $$
and so the MLE is defined as
$$\widehat{a}_{MLE} = \operatorname{arg} \operatorname{max}_{\boldsymbol{a}} \mathcal{L}(\boldsymbol{a} | n_{11},...,n_{1R},...n_{LR}, n_0).$$
Note that the log-likelihood function is not well-defined when the number of failed units due to cause $j,$ with $j=1,...,R,$ at any of the step-stress, is zero or $N.$ That is, for estimating the MLE, we need at least one failure due to each competing risk and under each of the stress levels.
For exact inference of the MLE conditional on that there is at least one failure observed due to each competing risk and under each of the stress levels, we refer to \cite{balakrishnan2007point}.
However, if the MLE is well defined, then it will satisfy the asymptotic properties, and particularly the asymptotic distribution, of the (unconditional) MLE. Also, the probability of not failing during a step-stress test tends to zero for each competing risk, and so asymptotically we will always have at least one failure per step and at each competing risk.
In what follows, we will assume that there is at least one observation for each combination of stress level and cause of failure. 

The MLE can also be justified using the information theory approach through the application of the Kullback-Leibler (KL) divergence. Divergence measures are designed to quantify the statistical similarity or dissimilarity between two distributions. Therefore, estimators based on divergences would be defined as the minimizers of a specific divergence measure between the theoretical distribution and the empirical distribution obtained from the data.
Let us consider the empirical estimator of the probability vector underlying, $\widehat{\boldsymbol{p}} = (p_{lj})_{l = 1,...,L, j=1,...,R} $ with $p_{lj} = \frac{n_{lj}}{N}.$ Then, the KL divergence between $\widehat{\boldsymbol{p}}$ and $\boldsymbol{p}(\boldsymbol{a})$ is given by
\begin{equation} \label{eq:KLdivergence}
	d_{KL}(\widehat{\boldsymbol{p}}, \boldsymbol{p}(\boldsymbol{a})) =  \sum_{j=1}^R\sum_{l=1}^{L}\widehat{p}_{lj} \log\left(\frac{\widehat{p}_{lj}}{p_{lj}(\boldsymbol{a})}\right) + \widehat{p}_{0} \log\left(\frac{\widehat{p}_{0}}{p_{0}(\boldsymbol{a})}\right)
\end{equation}
and consequently, the minimum KL divergence estimator is defined as
$$\widehat{a}_{KL} = \operatorname{arg} \operatorname{min}_{\boldsymbol{a}} d_{KL}(\widehat{\boldsymbol{p}}, \boldsymbol{p}(\boldsymbol{a})).$$
Removing all terms which do not depend on the parameter of the KL divergence in (\ref{eq:KLdivergence}), we readily obtain that the objective function defined by the KL divergence coincides with the negative log-likelihood of the model, and consequently the MLE coincides with the minimum KL divergence estimator. 

Building upon the idea  of using divergence-based estimators, and to overcome the robustness problem of the MLEs, we propose the use of the DPD for computing estimates of the model parameters.
The DPD between the empirical and theoretical probability vectors, $\widehat{\boldsymbol{p}}$ and $\boldsymbol{p}(\boldsymbol{a}),$ is given by,  for $\beta >0,$ 
\begin{equation}\label{eq:DPDloss}
	\begin{aligned}
	d_\beta(\widehat{\boldsymbol{p}}, \boldsymbol{p}(\boldsymbol{a}) ) = & \sum_{l=1}^L\sum_{j=1}^{R} \left(p_{lj}(\boldsymbol{a})^{1+\beta} -\left( 1+\frac{1}{\beta}\right) \widehat{p}_{lj}p_{lj}(\boldsymbol{a})^{\beta}  +\frac{1}{\beta} \widehat{p}_{lj}^{\beta+1} \right)\\
	&+ p_{0}(\boldsymbol{a})^{1+\beta} -\left( 1+\frac{1}{\beta}\right) \widehat{p}_{0}p_{0}(\boldsymbol{a})^{\beta}  +\frac{1}{\beta} \widehat{p}_{0}^{\beta+1},
\end{aligned}
\end{equation}
and the MDPDE is defined as
\begin{equation}\label{eq:MDPDE}
	\widehat{a}_{\beta} = \operatorname{arg} \operatorname{min}_{\boldsymbol{a}} d_\beta(\widehat{\boldsymbol{p}}, \boldsymbol{p}(\boldsymbol{a}) ) .
\end{equation}
Here, the tuning parameter $\beta$  controls the trade-off between efficiency and robustness; the lesser $\beta$, the more efficient the estimator will be, but less robust.
Indeed, the DPD can be defined at $\beta=0$  by taking continuous limits obtaining the KL divergence and thus, the MLE is the most efficient estimator, but least robust in the MDPDE family.

An important observation is that, since the DPD loss in Equation (\ref{eq:DPDloss}) is differentiable in $\boldsymbol{a}$ and the MDPDE is computed as its minimizer, the MDPDE must annul the first derivatives of the loss in (\ref{eq:DPDloss}). Therefore, the MDPDE estimating equations are given by
\begin{equation}\label{eq:estimating}
	\boldsymbol{W}^T\boldsymbol{D}_{\boldsymbol{p}(\boldsymbol{a})}^{\beta-1}\left( \widehat{\boldsymbol{p}}- \boldsymbol{p}(\boldsymbol{a})\right) = \boldsymbol{0}_{2R},
\end{equation}
where $\boldsymbol{0}_{2R}$ is the $2R$-dimensional null vector, $\boldsymbol{D}_{\boldsymbol{p}(\boldsymbol{a})}$ denotes a $(RL+1)\times(RL+1)$ diagonal matrix with diagonal entries $p_{lj}(\boldsymbol{a}),$ $l=1,...,L, j=1,...,R,$ and $p_0(\boldsymbol{a}),$ and  $\boldsymbol{W}$ is a $(L R+1) \times 2 R$ matrix 
consisting of $R$  columns blocks of dimension $(LR+1)\times 2$ defined as
$\boldsymbol{w}^{k} = \left[ (\boldsymbol{z}_{l-1,j}^{k}-\boldsymbol{z}_{lj}^{k})^T_{j=1,..R, l=1,...,L}, \boldsymbol{z}_0^{k}\right],$ for $k=1,...,R$ 
with $\boldsymbol{z}_{l,j}^{k} = ((\boldsymbol{z}_{l,j})_{1k}, (\boldsymbol{z}_{lj})_{2k}),$
\begin{equation*}
	\begin{aligned}
		(\boldsymbol{z}_{lj})_{1k} &= 
		\begin{cases}
			\pi_{ij} \exp\left(-\sum_{j=1}^{R} \frac{IT_{l}+h_j^{(i)}}{\exp(a_{0j}+a_{1j}x_i)}\right) \left(\pi_{ik}+\frac{IT_l + h_k^{(i)}}{\exp(a_{0k}+a_{1k}x_i)}\right), & j \neq k, \\
			\pi_{ij}\exp\left(-\sum_{j=1}^{R} \frac{IT_{l}+h_j^{(i)}}{\exp(a_{0j}+a_{1j}x_i)}\right)\left(\pi_{ij}-1+\frac{IT_l + h_j^{(i)}}{\exp(a_{0j}+a_{1j}x_i)}\right),  & j = k,\\
		\end{cases}\\
		(\boldsymbol{z}_{lj})_{2k} &= 
		\begin{cases}
			\pi_{ij}\exp\left(-\sum_{j=1}^{R} \frac{IT_{l}+h_j^{(i)}}{\exp(a_{0j}+a_{1j}x_i)}\right)
			\left(\pi_{ik}x_i+\frac{-h_k^{\ast(i)} + (IT_l + h_k^{(i)})x_i}{\exp(a_{0k}+a_{1k}x_i)}\right), & j \neq k, \\
			\pi_{ij}\exp\left(-\sum_{j=1}^{R} \frac{IT_{l}+h_j^{(i)}}{\exp(a_{0j}+a_{1j}x_i)}\right)
			\left(\pi_{ij}x_i-x_i+\frac{-h_j^{\ast(i)} + (IT_l + h_j^{(i)})x_i}{\exp(a_{0j}+a_{1j}x_i)}\right),  & j = k,\\
		\end{cases}
	\end{aligned}
\end{equation*}
where $i$ is chosen such that $\tau_{i-1} \leq IT_{l} <\tau_i,$ $h_j^{ (i)}$ is as defined in (\ref{eq:shiftingtime2}) and
\begin{equation}
	h_j^{\ast (1)} = 0 \hspace{0.3cm} \text{ and } \hspace{0.3cm} h_j^{\ast (2)} = \tau_1\frac{\exp(a_{0j}+a_{1j}x_2)}{\exp(a_{0j}+a_{1j}x_1)} (x_2-x_1).
\end{equation}
For the last column of the block $\boldsymbol{w}^{k},$ the entries of the matrix are defined as
$$(\boldsymbol{z}_0)_{1k} = \exp\left(- \sum_{j=1}^R \frac{\tau_2+h_j^{(2)}}{\exp(a_{0j} + a_{1j}x_2)} \right) \left(\frac{\tau_2+h_j^{(2)}}{\exp(a_{0k}+a_{1k}x_2)}\right),$$

$$(\boldsymbol{z}_0)_{2k} = \exp\left(- \sum_{j=1}^R \frac{\tau_2+h_j^{(2)}}{\exp(a_{0j} + a_{1j}x_2)} \right) \left(\frac{-h_{k}^{\ast (2)} + (\tau_2 + h_k^{(2)})x_2}{\exp(a_{0k}+a_{1k}x_2)}\right).$$
For detailed derivations of the above matrices, see Appendix \ref{AppA}.

The above equations characterize the MDPDE family as M-estimators, and so they would enjoy all desirable properties of the M-estimators such as consistency and asymptotic normality. Further, the asymptotic distribution of the MDPDE in general statistical models was established in \cite{basu1998robust}. The next result states the asymptotic distribution for the interval-monitoring step-stress ALT model with competing risks and exponential lifetimes.

\begin{theorem}\label{thm:asymptotic}
	Let $\boldsymbol{a}_0$ be the true value of the interval-monitoring step-stress model parameter $\boldsymbol{a}$ and consider the MDPDE of $\boldsymbol{a}$ with tuning parameter $\beta$, $\boldsymbol{\widehat{a}}^{\beta}.$ 
	Then, the asymptotic distribution of $\boldsymbol{\widehat{a}}^{\beta}$ is given by
	$$ \sqrt{N}\left(\boldsymbol{\widehat{a}}^{\beta} - \boldsymbol{a}_0\right) \xrightarrow[N \rightarrow \infty]{L}\mathcal{N}\left(\boldsymbol{0}, \boldsymbol{J}_\beta^{-1}(\boldsymbol{a}_0)\boldsymbol{K}_\beta(\boldsymbol{a}_0)\boldsymbol{J}_\beta^{-1}(\boldsymbol{a}_0)\right),$$
	where
	\begin{equation}\label{eq:JK}
		\boldsymbol{J}_\beta(\boldsymbol{a}_0) = \boldsymbol{W}^T \boldsymbol{D}_{\boldsymbol{p}(\boldsymbol{a_0})}^{\beta-1} \boldsymbol{W},
		\hspace{0.3cm}  \hspace{0.3cm}
		\boldsymbol{K}_\beta(\boldsymbol{a}_0) = \boldsymbol{W}^T \left( \boldsymbol{D}_{\boldsymbol{p}(\boldsymbol{a_0})}^{2\beta-1}-\boldsymbol{p}(\boldsymbol{a}_0)^{\beta}\boldsymbol{p}(\boldsymbol{a}_0)^{\beta T}\right) \boldsymbol{W},
	\end{equation}
	with $D_{\boldsymbol{p}(\boldsymbol{a_0})}$ denoting the diagonal matrix with entries $p_{lj}(\boldsymbol{a_0}),$ $l=1,...,L, j=1,...,R,$ and $\boldsymbol{p}(\boldsymbol{a}_0)^{\beta}$ denoting the vector with components $p_{lj}(\boldsymbol{a}_0)^{\beta}.$
\end{theorem}

Result \ref{thm:asymptotic} can be used for finding approximate standard errors of the MDPDEs as well as approximate confidence intervals for the model parameters. 
Indeed,  the diagonal entries of the covariance matrix 
\begin{equation} \label{eq:Sigma}
	\boldsymbol{\Sigma}(\boldsymbol{a}_0) = \boldsymbol{J}_\beta^{-1}(\boldsymbol{a}_0)\boldsymbol{K}_\beta(\boldsymbol{a}_0)\boldsymbol{J}_\beta^{-1}(\boldsymbol{a}_0)
\end{equation} are the approximate standard errors of the estimators. 
Besides, since the MDPDEs are consistent estimators of $\boldsymbol{a}_0,$ we can readily obtain consistent estimators of the error by plugging-in the MDPDEs in $\boldsymbol{\Sigma}(\boldsymbol{a})$ and the  $100(1-\alpha)\%$ approximate confidence intervals for $a_{ij}$ can be obtained as
\begin{equation}\label{asymptotic}
	IC_\alpha(a_{ij}) = \left[ \widehat{a}^\beta_{ij} \pm \frac{z_{\alpha/2}}{N} \Sigma(\widehat{\boldsymbol{a}}^\beta)_{ij}\right], \hspace{0.3cm} i = 0,1 \hspace{0.1cm} j = 1,...,R,
\end{equation}
where $z_{\alpha}$ denotes the upper $\alpha$ quantile of the standard normal distribution and $\boldsymbol{\Sigma}(\boldsymbol{a})$ is as defined in Equation (\ref{eq:Sigma}). 



\section{Robust point estimation and confidence intervals for lifetime characteristics \label{sec:char}}

In many reliability analyses, one is interested in inferring a lifetime characteristic rather than the complete lifetime distribution function.
The lifetime characteristics may  help in understanding the reliability characteristics of products and guides in decision-making processes related to maintenance, replacement, and risk management. In this section, we develop point estimation and approximate confidence intervals for the mean lifetime to failure (MTTF) under constant stress, distributional quantiles, and reliability of a product at a mission time $t_0.$ For a more accurate performance of the confidence interval, we also develop transformed confidence intervals for the aforementioned lifetime characteristics. Additionally, we present a bootstrap algorithm for building confidence intervals as an alternative to asymptotic intervals.

Let us denote $x_0$ for the stress level under normal operating conditions. For partially accelerated step-stress plans, $x_0$ coincides with the first stress level considered, denoted by $x_1$ above. 
Given MDPDE of the step-stress model parameters, $\widehat{\boldsymbol{a}}^\beta,$ we can readily estimate the PDF of the lifetime to failure by plugging-in the estimators in Equation (\ref{eq:PDFmarginal}), yielding 
%
\begin{equation}
	f_T(t) 
	=
	\left[\sum_{j=1}^R \frac{1}{\widehat{\theta}^\beta_{0j}}\right] \exp\left(- \sum_{j=1}^R \frac{t}{\widehat{\theta}^\beta_{0j}}\right), 
\end{equation}
where the exponential parameters are estimated as $\widehat{\theta}^\beta_{0j} = \widehat{a}_{0j} + \widehat{a}_{1j}x_0$ for each competing risk $j=1,..,R.$
The distribution of the lifetime to failure under constant stress $x_0,$ $T,$ is an exponential distribution with failure rate $1/\theta = 
\sum_{j=1}^R 1/\widehat{\theta}^\beta_{0j}$
and thus using the above PDF, we can compute an estimate for any lifetime characteristic. Moreover, using the asymptotic standard error of the MDPDEs, we can also obtain approximate standard error for the lifetime characteristic of interest.

\subsection{Point estimation and confidence intervals for the mean lifetime to failure}

The mean lifetime to failure (MTTF) of the lifetime $T$ under normal operating conditions is given by
\begin{equation}\label{eq:MTTF}
		E(\boldsymbol{a}) = E[T] 
		= \left[\sum_{j=1}^R \exp(-a_{0j}-a_{1j}x_0)\right]^{-1}
\end{equation}
and consequently, the MDPDE of the MTTF under normal operating conditions can be obtained as 
\begin{equation}
	E(\widehat{\boldsymbol{a}}^\beta)
	= \left[\sum_{j=1}^R \exp(-\widehat{a}^\beta_{0j}-\widehat{a}^\beta_{1j}x_0)\right]^{-1}.
\end{equation}

Using the asymptotic distribution of the MDPDE stated in Result \ref{thm:asymptotic} and applying the delta method in Equation (\ref{eq:MTTF}), we have the asymptotic distribution of the MDPDE for the MTTF as

$$\sqrt{N}\left(	E(\widehat{\boldsymbol{a}}^\beta) - 	E(\boldsymbol{a}) \right) \xrightarrow[N\rightarrow \infty]{L} \mathcal{N}\left(0, \sigma^2(E) \right)$$ 
with 
$$\sigma^2(E) = \nabla E(\boldsymbol{a}_0)^T \boldsymbol{\Sigma}_\beta(\boldsymbol{a}_0) \nabla E(\boldsymbol{a}_0),$$
where $\nabla E(\boldsymbol{a}_0)$ is the $2R$-dimensional gradient of the function $E(\cdot)$ whose components are given by $2-$dimensional blocks of the form
$$\nabla E(\boldsymbol{a}_0)_{2k-1,2k} = \left(\frac{\pi_{0k}}{\sum_{j=1}^R \exp(-a^\beta_{0j}-a^\beta_{1j}x_0)} \left(1, x_0\right)^T \right)_{k=1,...,R},$$
and $\boldsymbol{\Sigma}_\beta(\boldsymbol{a})$ is the asymptotic covariance matrix of the MDPDE given in (\ref{eq:Sigma}).
%
%
Further, since the MDPDE is a consistent estimator of the true value of the parameter $\boldsymbol{a}_0$, we can estimate the standard error of the MDPDE for the MTTF as 
$$\widehat{\sigma}^2(E) = \nabla E(\widehat{\boldsymbol{a}}^\beta)^T \boldsymbol{\Sigma}_\beta(\widehat{\boldsymbol{a}}^\beta) \nabla E(\widehat{\boldsymbol{a}}^\beta)$$
and so $100(1-\alpha)\%$ approximate confidence interval is given by 
$$IC_\alpha(E) = \left[ E(\widehat{\boldsymbol{a}}^\beta) \pm z_{\alpha/2} \frac{\widehat{\sigma}(E)}{N} \right].$$
Note that the MTTF must be positive and therefore the above interval should be truncated for negative lower bounds.

\subsection{Point estimation and confidence intervals for the reliability of the devices}

Another lifetime characteristic of great interest for manufacturers is product reliability at a certain mission time.
Many industries, such as aerospace, defense, and critical infrastructure, have specific reliability requirements for systems or components during mission-critical operations. In those cases, an estimate of the reliability is necessary to ensure that the requirements of the product for sale are met.
Moreover, it helps the manufacturer to ensure that his/her product will perform satisfactorily during the intended mission or operation, and 
consequently to determine an appropriate time for warranty coverages with a low failure rate.
%
%
As aforementioned, a robust estimate of the reliability function of the device lifetime $T$ and marginal lifetime due to risk $j,$ $T_j,$ with $j=1,...,R,$  can be directly determined by plugging-in the MDPDEs into their theoretical expressions. Here, we derive point estimation of the reliability of the device and reliability of each component at mission time $t_0$ as well as the corresponding approximate confidence intervals. 

Let us first consider the estimated reliability of the lifetime $T$ under the normal operating conditions $x_0$ at a mission time $t_0$ given by
\begin{equation}\label{eq:R0}
	R_0(\boldsymbol{a}) 
	= \exp\left(- t_0\sum_{j=1}^R \exp(-a_{0j}-a_{1j}x_0)\right), 
\end{equation}
and its corresponding MDPDE given by
\begin{equation}
	R_0(\widehat{\boldsymbol{a}}^\beta) 
	= \exp\left(- t_0\sum_{j=1}^R \exp(-\widehat{a}^\beta_{0j}-\widehat{a}^\beta_{1j}x_0)\right). 
\end{equation}
Using the asymptotic distribution of the MDPDE stated in Result \ref{thm:asymptotic} and applying the delta method in Equation (\ref{eq:R0}), we have the asymptotic distribution of the MDPDE for the reliability of the device at mission time $t_0$ as

$$\sqrt{N}\left(R_0(\widehat{\boldsymbol{a}}^\beta) - 	R_0(\boldsymbol{a}) \right) \xrightarrow[N\rightarrow \infty]{L} \mathcal{N}\left(0, \sigma^2(R_0) \right),$$ 
with 
$$\sigma^2(R_0) = \nabla R_0(\boldsymbol{a}_0)^T \boldsymbol{\Sigma}_\beta(\boldsymbol{a}_0) \nabla R_0(\boldsymbol{a}_0),$$
where $\nabla R_0(\boldsymbol{a})$ the $2R$-dimensional gradient of the reliability function at mission $t_0$ with components formed by $2$-dimensional vectors
$$\nabla R_0(\boldsymbol{a})_{2k-1,2k} = R_0(\boldsymbol{a})t_0\left(\exp(-a_{0k}-a_{1k}x_0), \exp(-a_{0k}-a_{1k}x_0)x_0\right)^T, $$ and $\boldsymbol{\Sigma}_\beta(\boldsymbol{a})$ is as defined in Equation (\ref{eq:Sigma})
Finally, a robust and consistent estimate of the standard errors of $\sigma^2(R_0)$  can be easily obtained by using MDPDEs and thus approximate $100(1-\alpha)\%$ confidence intervals for the reliability of the product as
$$IC_\alpha(R_0) = \left[ R_0(\widehat{\boldsymbol{a}}^\beta) \pm z_{\alpha/2} \frac{\widehat{\sigma}(R_0)}{N} \right],$$ 
where $z_{\alpha}$ denotes the upper $\alpha$ quantile of a standard normal distribution. Because the reliability of a product is bounded by 0 (lower bound) and 1 (upper bound), the above confidence interval may have to be truncated.


\subsection{Point estimation and confidence intervals for the distribution quantiles}

Lifetime quantiles can also be used to set reliability targets or specifications for a product.
By specifying a desired quantile value, manufacturers can establish performance requirements and ensure that their products meet certain reliability standards.
Quantile estimation can serve for inferring central tendencies such as the median of the lifetime distribution, which is often the target of a reliability analysis, but can also can serve for estimating limit times at the distribution tail. Moreover, estimating those tail quantiles may help in  identifying any abnormal performance of new devices.

 Let us consider the lower $\alpha_0$ quantile of the lifetime distribution $T$ given by
$$Q_{1-\alpha_0}(\boldsymbol{a}) = \frac{-\log(1-\alpha_0)}{\sum_{j=1}^R \exp(-a_{0j}-a_{1j}x_0)}$$
and its corresponding MDPDE given by $$Q_{1-\alpha_0}(\widehat{\boldsymbol{a}}^\beta)= \frac{-\log(1-\alpha_0)}{\sum_{j=1}^R \exp(-\widehat{a}^\beta_{0j}-\widehat{a}^\beta_{1j}x_0)}.$$

We can derive the asymptotic distribution of the above estimator by applying again the delta method yielding the asymptotic result as
$$\sqrt{N}\left(Q_{1-\alpha_0}(\widehat{\boldsymbol{a}}^\beta) - 	Q_{1-\alpha}(\boldsymbol{a}) \right) \xrightarrow[N\rightarrow \infty]{L} \mathcal{N}\left(0, \sigma^2(Q_{1-\alpha}) \right)$$ 
where 
\begin{equation*}
	\sigma^2(Q_{1-\alpha})  =  \nabla Q_{1-\alpha_0}(\boldsymbol{a}_0)^T \boldsymbol{\Sigma}_\beta(\boldsymbol{a}_0) \nabla Q_{1-\alpha}(\boldsymbol{a}_0),
\end{equation*}
and $\nabla Q_{1-\alpha_0}(\boldsymbol{a}_0)$ is the $2R$-dimensional gradient of the function $Q_{1-\alpha_0}(\cdot)$ whose components are given by $2-$dimensional blocks 
$$\nabla Q_{1-\alpha_0}(\boldsymbol{a}_0)_{2k-1,2k} =  \left(\frac{-\log(1-\alpha_0)\pi_{0k}}{\sum_{j=1}^R \exp(-a^\beta_{0j}-a^\beta_{1j}x_0)} \left(1, x_0\right)^T \right)_{k=1,...,R},$$
and $\boldsymbol{\Sigma}_\beta(\boldsymbol{a})$ is the asymptotic covariance matrix of the MDPDE given in (\ref{eq:Sigma}).
Then, with the previous notation, an approximate $100(1-\alpha)\%$ confidence interval for an $1-\alpha_0$ quantile is then given by
$$IC_\alpha(Q_{1-\alpha_0}) = \left[ Q_{1-\alpha_0}(\widehat{\boldsymbol{a}}^\beta) \pm z_{\alpha/2} \frac{\widehat{\sigma}(Q_{1-\alpha_0})}{N} \right].$$ 
Note that $\nabla Q_{1-\alpha_0}(\boldsymbol{a}_0)_{2k-1,2k} = -\log(1-\alpha_0) \nabla E(\boldsymbol{a}_0)_{2k-1,2k} $ and so the asymptotic standard  error of the quantile coincides with the asymptotic standard  error of the mean lifetime up to a factor of $\log(1-\alpha_0).$
\subsection{Transformed confidence intervals \label{sec:transformedIC}}

The above asymptotic confidence intervals are based on the asymptotic properties of the MDPDEs and so they may be satisfactory only for large sample sizes.
As we have mentioned above, the lifetime characteristics must satisfy some natural restrictions and therefore the confidence intervals may have to be truncated.
In this regard, \cite{viveros1993statistical} employed a logarithmic transformation to the mean lifetime and distribution quantiles and a logit transformation to the reliability to obtain more accurate confidence intervals based on the MLEs. These transformations were also applied in \cite{balakrishnan2022robust} to the general MDPDEs for the step-stress ALT, producing the same effect for the robust confidence intervals. We also adopt the transformation approach here for improving the coverage of the approximate intervals in small samples.

Let us consider the following transformed characteristics;
\begin{equation*}\label{eq:transformed}
		\phi_E(\boldsymbol{a}) = \log\left(E(\boldsymbol{a})\right), \hspace{0.3cm}
		\phi_R(\boldsymbol{a}) = \operatorname{logit}\left(R_0(\boldsymbol{a})\right) 
		\hspace{0.3cm} \text{and} \hspace{0.3cm}
		\phi_Q(\boldsymbol{a}) = \log\left(Q_{1-\alpha_0}(\boldsymbol{a})\right).
\end{equation*}
Applying the delta method to the asymptotic distributions of the mean lifetime, reliability, and distribution quantile, we can obtain confidence intervals for $	\phi_E(\boldsymbol{a}), \phi_R(\boldsymbol{a}),$ and $\phi_Q(\boldsymbol{a}).$ 
Inverting the injective transformations, transformed confidence intervals for the mean lifetime, reliability at a mission time $t_0$ and $1-\alpha_0$ upper quantile can be obtained, respectively, as
\begin{equation}\label{eq:transformedci}
	\begin{aligned}
	IC_\alpha(E) = 	&\left[E(\widehat{\boldsymbol{a}}^\beta) \exp\left(-\frac{z_{\alpha/2}}{\sqrt{N}} \frac{\sigma(E)}{E(\widehat{\boldsymbol{a}}^\beta)}\right), 
		E(\widehat{\boldsymbol{a}}^\beta) \exp\left(\frac{z_{\alpha/2}}{\sqrt{N}} \frac{\sigma(E)}{E(\widehat{\boldsymbol{a}}^\beta)}\right)\right],\\
	IC_\alpha(R_0) = 	&\left[\frac{R_0(\widehat{\boldsymbol{a}}^\beta)}{R_0(\widehat{\boldsymbol{a}}^\beta) + \left(1-R_0(\widehat{\boldsymbol{a}}^\beta)\right)S},
		\frac{R_0(\widehat{\boldsymbol{a}}^\beta)}{R_0(\widehat{\boldsymbol{a}}^\beta) + \left(1-R_0(\widehat{\boldsymbol{a}}^\beta)\right)/S}
		\right],\\
	IC_\alpha(Q_{1-\alpha_0}) = 	&\left[Q_{1-\alpha_0}(\widehat{\boldsymbol{a}}^\beta) \exp\left(-\frac{z_{\alpha/2}}{\sqrt{N}} \frac{\sigma(Q_{1-\alpha_0})}{Q_{1-\alpha_0}(\widehat{\boldsymbol{a}}^\beta)}\right), 
		Q_{1-\alpha_0}(\widehat{\boldsymbol{a}}^\beta) \exp\left(\frac{z_{\alpha/2}}{\sqrt{N}} \frac{\sigma(Q_{1-\alpha_0})}{Q_{1-\alpha_0}(\widehat{\boldsymbol{a}}^\beta)}\right)\right],
	\end{aligned}
\end{equation}
where 
$ S= \exp \left(\frac{z_{\alpha/2}}{\sqrt{N}} \frac{\sigma(R_0)}{R_0(\widehat{\boldsymbol{a}}^\beta) (1-R_0(\widehat{\boldsymbol{a}}^\beta))}\right),$ $z_{\alpha}$ denotes the upper $\alpha$ quantile of a standard normal distribution and $\sigma(E), \sigma(R_0)$ and $\sigma(Q_{1-\alpha})$ are the estimated standard errors of the mean lifetime, reliability at mission time $t_0$ and $1-\alpha_0$ quantile, respectively.

\subsection{Bootstrap confidence intervals}

An appealing parametric alternative in case of small sample sizes to asymptotic confidence intervals are the so-called bias-corrected and accelerated (BCa) percentile bootstrap confidence intervals. 


For each of the competing risks under study, a set of $N$ failure times are simulated from the estimated marginal distributions. Since only the smallest lifetime is observed in practice, the observed lifetime is computed as the minimum lifetime on the marginal competing risks  $T=\operatorname{min}(T_1,...,T_R)$ and the associated cause of failure  is recorded. 
The specific algorithm is as follows:

\vspace{0.3cm}
\noindent \textbf{Algorithm 1: Bootstrap confidence interval}\\
\hrule

\begin{enumerate}
	\item  
	 Obtain the MDPDEs defined in (\ref{eq:DPDloss})  based on the observed simple step-stress  interval-censored sample of size $N$;
	
	\item  Generate $R$ random samples from independent and identically distributed uniform distributions, $\mathcal{U}(0, 1),$ of size $N.$ We denote  $\boldsymbol{U}^{(j)} = \{U^{(j)}_1,...,U_N^{(j)}\}$ the $j$-th random sample;
	
	\item Transform each uniform sample $\boldsymbol{U}^{(j)}$ into a lifetime vector $\boldsymbol{T}^{(j)} = (T_1^{(j)},...,T_N^{(j)})$ by inverting the marginal CDFs of the competing risk as follows:
	\begin{equation*}
		\begin{cases}
			T_{k}^{(j)} = -\log(1-U_{k}^{(j)})\widehat{\theta}_{1j}^{1}, & U_k^{(j)} < F_j(\tau_1),\\
			T_{k}^{(j)} = -\log(1-U_{k}^{(j)})\widehat{\theta}_{2j}^{1} + h_j^{(2)}, &  F_j(\tau_1) \leq U_k^{(j)},\\
		\end{cases}
	\end{equation*}
	with $k=1,...,N;$
	\item For each sample index $k,$ 
	compute the simulated failure time as the minimum of the marginal failure times, $T_k = \operatorname{min}(T_k^{(1)},...,T_k^{(R)}),$ and record the index of the minimum, denoted by $C_{T_k};$
	
	\item Obtain the order statistics of the lifetime sample with their associated cause of failure
	
	$$ \{(T_{1:N}, C_{T_{1:N}}), ... ,(T_{N:N}, C_{T_{N:N}})\};$$
	
	\item Based on the above simulated sample, compute the number of  failures in each inspection interval, i.e., the frequencies $n_{ij}, i=1,2,$  $j=1,...,R,$ as well as the number of surviving units $n_0.$  If the number of failures due to a competing risk at any of the step-stress is zero on a certain simulated dataset, discard this dataset and generate a new one;
	
	\item  Calculate the new MDPDEs of the model parameters $\boldsymbol{a}$ based on the simulated sample from (\ref{eq:MDPDE}), denoted by $(\widehat{\boldsymbol{a}}^\beta)_{b}.$ Using the new MDPDEs, compute the lifetime characteristic of interest, either the mean lifetime, reliability, or distribution quantile. We denote by $\widehat{S}^\beta_{b}$ the estimated value of the characteristic and $S$ its true value;
	
	\item Repeat Steps 2 through 7 $B$ times. 
	Then, arrange all the values of $\widehat{S}^\beta_{b}$ in an ascending order to obtain the bootstrap sample
	$$\{\widehat{S}^\beta_{1:B} < \cdots < \widehat{S}^\beta_{B:B}\}.$$
\end{enumerate}	
Using the bootstrap sample generated by the algorithm given above, the two-sided $100(1-\alpha)\%$ BCa percentile bootstrap CI  is obtained as
$$IC_{1-\alpha}(S) = \left[\widehat{S}^\beta_{[\gamma_1B]:B}, \widehat{S}^\beta_{[\gamma_2B]:B} \right],$$
where
$$\gamma_1 = \Phi\left(\widehat{z}_{0} + \frac{\widehat{z}_0-z_{\alpha/2}}{1-\widehat{\gamma}(\widehat{z}_0-z_{\alpha/2})}\right)
\hspace{0.3cm} \text{and} \hspace{0.3cm}
\gamma_2 = \Phi\left(\widehat{z}_{0} + \frac{\widehat{z}_0+z_{\alpha/2}}{1-\widehat{\gamma}(\widehat{z}_0+z_{\alpha/2})}\right),$$
 and $\Phi(\cdot)$ is the standard normal CDF with upper $\alpha$ quantile $z_{\alpha}.$ Furthermore, $\widehat{z}_0$ and $\widehat{\gamma}$  stand for the bias correction and the acceleration, the first is estimated by
$$\widehat{z}_0 = \Phi^{-1}\left(\frac{\# \text{ of } \{(\widehat{S}_{b}^\beta) \leq \widehat{S}^\beta\} }{B}\right)$$ 
 and according to \cite{effron1993introduction}, a suggested estimate of the acceleration bias is given by
$$
\widehat{\gamma} = \frac{1}{6}\left[\sum_{l=1}^{N_L} \left( \widehat{S}_{l}^\beta - \widehat{S}_{(\cdot)}^\beta \right)^2\right]^{-3/2}  \left[\sum_{l=1}^L \left( \widehat{S}_{l} - \widehat{S}_{(\cdot)}^\beta \right)^3\right],
$$
where $\widehat{S}_{l}^\beta$
 is the MDPDE of the lifetime characteristic $S$ with tuning parameter $\beta$  based on the initial observed sample, but with the $l$-th observation deleted (i.e. the jackknife estimate), and
 $$\widehat{S}_{(\cdot)}^\beta = \frac{1}{N_L}\sum_{l=1}^{N_L} \widehat{S}_{l}^\beta$$
with $N_L$ being the total number of failures observed.

The above bootstrap algorithm can be used for building BCa confidence intervals of the model parameters, $a_{ij}$.

\section{Cause-specific lifetime characteristics \label{sec:specificchar}}

It can be of interest to separately estimate the lifetime characteristics of a specific failure cause. This is particularly relevant in clinical analyses, wherein the causes of death can be attributed to different diseases, and analyzing the distribution characteristics of each disease separately can provide valuable insight. Besides, estimating lifetime patterns of the different competing risks individually has also extensive applications in industrial experiments, where identifying the unique patterns and risks associated with each risk aids in making informed decisions regarding the replacement or maintenance of specific components, optimizing reliability, and ensuring efficient operation of industrial systems.

Let us consider the marginal PDF of the lifetime to failure due to cause $j,$ for a fixed $j,$ given in Equation (\ref{eq:PDFmarginal}).  Given the MDPDEs related to cause $j,$ $\widehat{\boldsymbol{a}}_j^\beta = (\widehat{a}_{0j}^\beta, \widehat{a}_{1j}^\beta),$ the marginal PDF at constant stress $x_i$ can be readily estimated as 
\begin{equation}\label{eq:estimatedPDFmarginal}
	f_j(t) = \frac{1}{\widehat{\theta}_{ij}^\beta}\exp \left(-\frac{t_j}{\widehat{\theta}^\beta_{ij}} \right),
\end{equation} 
where $\widehat{\theta}^\beta_{ij} = \widehat{a}_{0j}^\beta + \widehat{a}_{0j}^\beta x_{i}.$ From the above exponential PDF, we can compute any lifetime characteristic such as mean, median, reliability, or distribution quantiles. The corresponding estimates will inherit the properties of the MDPDEs in terms of efficiency and robustness. Thus, estimates of the lifetime characteristics using small values of $\beta$ will be more efficient but less robust, and in contrast, large values of $\beta$ would provide more robust  but less efficient estimators.

\subsection{Point estimation and confidence intervals for cause-specific mean lifetime}

Marginal lifetime follows exponential distribution and so the estimated scale parameter of the exponential distribution for a fixed stress, say $x_0,$ determines its mean lifetime to failure (MTTF) as
 \begin{equation}
 	E_j(\widehat{\boldsymbol{a}}_j^\beta) = \widehat{\theta}^\beta_{0j} = \exp(\widehat{a}_{0j}^\beta + \widehat{a}_{0j}^\beta x_{0}).
 \end{equation} 
Again, using the delta method and the asymptotic distribution of the MDPDE $\widehat{\boldsymbol{a}}_j^\beta,$ we have 

$$\sqrt{N}\left(E_j(\widehat{\boldsymbol{a}}_j^\beta)- E_j(\boldsymbol{a}_j^\beta)\right)\xrightarrow[N\rightarrow \infty]{L} \mathcal{N}\left(0, \sigma^2(E_j) \right)$$  
with 
$$\sigma^2(E_j) = E_j(\widehat{\boldsymbol{a}}_j^\beta)^2\left(1,x_0\right)^T \left[\boldsymbol{\Sigma}_\beta(\boldsymbol{a}_0)\right]_{2j-1,j} \left(1,x_0\right),$$
where 
$\left[\boldsymbol{\Sigma}_\beta(\boldsymbol{a})\right]_{2j-1,j}$ is the $2\times 2$ submatrix
of the matrix $\boldsymbol{\Sigma}_\beta(\boldsymbol{a})$ as defined  in Equation (\ref{eq:Sigma}) at entries $(2j-1,j)\times (2j-1,j).$

\subsection{Point estimation and confidence intervals for cause-specific reliability}

In some industrial products, such as electronic devices, the competing risks are different components of the devices which may fail and it would be also of interest to estimate the marginal reliability of these components causing failure.  In those cases, understanding the marginal reliability of each component may help in making decisions about future product replacements.
We consider the marginal reliability of  a fixed  competing risk $j$, with $j=1,...,R,$ under normal operating conditions defined by
$$R_j(\boldsymbol{a}_j) = \exp \left(-t_j\exp(-a_{0j}-a_{1j}x_0) \right).
$$

Following arguments similar to those for the lifetime of the device, the MDPDE of the reliability of the $j-$th component is given by $R_j(\widehat{\boldsymbol{a}}_j)$ and it asymptotic distribution is 
$$\sqrt{N}\left(R_j(\widehat{\boldsymbol{a}}_j^\beta) - 	R_j(\boldsymbol{a}_j) \right) \xrightarrow[N\rightarrow \infty]{L} \mathcal{N}\left(0, \sigma^2(R_j) \right)$$ 
with 
$$\sigma^2(R_j) = R_j(\boldsymbol{a}_j)^2 t_0^2 \exp(-2a_{0j}-2a_{1j}x_0)\left(1,x_0\right)^T [\boldsymbol{\Sigma}_\beta(\boldsymbol{a})]_{2j-1,2j} \left(1,x_0\right),$$
where $[\boldsymbol{\Sigma}_\beta(\boldsymbol{a}_j)]_{2j-1,2j}$ is the  $2\times2$ diagonal block of the asymptotic covariance matrix given in Equation (\ref{eq:Sigma}) at entries $(2j-1,2j)\times(2j-1,2j).$
Accordingly, a $100(1-\alpha)$ approximate confidence interval for the marginal reliability of each competing risk is given by
$$ IC_\alpha(R_j) = \left[ R_j(\widehat{\boldsymbol{a}}_j^\beta) \pm z_{\alpha/2} \frac{\widehat{\sigma}(R_j)}{N} \right],$$
where $z_{\alpha}$ denotes the upper $\alpha$ quantile of a standard normal distribution. 
Note that the reliability of each competing risk lies in the interval $[0,1]$ and so the above confidence interval may be truncated to ensure that condition is fulfilled. Using transformed confidence intervals as in Section (\ref{sec:transformedIC}) we have the transformed confidence interval as 

\begin{equation*}
		IC_\alpha(R_j) = 	\left[\frac{R_j(\widehat{\boldsymbol{a}}_j^\beta)}{R_0(\widehat{\boldsymbol{a}}_j^\beta) + \left(1-R_j(\widehat{\boldsymbol{a}}_j^\beta)\right)S},
		\frac{R_j(\widehat{\boldsymbol{a}}_j^\beta)}{R_j(\widehat{\boldsymbol{a}}^\beta) + \left(1-R_j(\widehat{\boldsymbol{a}}_j^\beta)\right)/S}
		\right],\\
\end{equation*}
where 
$ S= \exp \left(\frac{z_{\alpha/2}}{\sqrt{N}} \frac{\sigma(R_j)}{R_j(\widehat{\boldsymbol{a}}_j^\beta) (1-R_j(\widehat{\boldsymbol{a}}_j^\beta))}\right).$

\subsection{Point estimation and confidence intervals for cause-specific quantiles}

Marginal distribution quantiles determine a limit time when a component of the device (the competing risk) is expected to survive with a fixed probability.
Therefore, estimating marginal tail quantiles is useful for determining times of change of certain components. Additionally, fixing a probability of $\alpha_0=1/2,$ we can estimate the median life of each component separately.

Because the marginal lifetime distribution under constant stress $x_0$ of any $j$-component follows an exponential distribution with scale parameter $\theta_{0j} = \exp(a_{0j}+a_{01}x_0)$, the estimated marginal lower $\alpha_0$ quantile is given by
\begin{equation}
	Q_{j,1-\alpha_0}(\widehat{\boldsymbol{a}}_{j}) = -\log(1-\alpha_0)\exp(\widehat{a}_{0j}+\widehat{a}_{01}x_0)
\end{equation}
and by the delta method, the asymptotic distribution of the quantile estimate based on the MDPDE, with tuning parameter $\beta,$ is given by
$$\sqrt{N}\left(Q_{j,1-\alpha_0}(\widehat{\boldsymbol{a}}_{j}) - Q_{j,1-\alpha_0}(\boldsymbol{a}_{j}) \right) \xrightarrow[N\rightarrow \infty]{L} \mathcal{N}\left(0, \sigma^2(Q_{j,1-\alpha_0}) \right),$$ 
with $\boldsymbol{a}_j$ being the true value of the parameter and 
$$\sigma^2(Q_{j,1-\alpha_0}) =  Q_{j,1-\alpha_0}(\boldsymbol{a}_{j})\left(1,x_0\right)^T [\boldsymbol{\Sigma}_\beta(\boldsymbol{a})]_{2j-1,2j} \left(1,x_0\right),$$
with $[\boldsymbol{\Sigma}_\beta(\boldsymbol{a}_j)]_{2j-1,2j}$ being a   $2\times2$ diagonal block of the asymptotic covariance matrix given in Equation (\ref{eq:Sigma}) at entries $(2j-1,2j)\times(2j-1,2j).$

The direct asymptotic and transformed $100(1-\alpha)\%$ asymptotic confidence intervals for the cause-specific $\alpha_0$ quantile are given by
$$ IC_\alpha(Q_{j,1-\alpha_0}) = \left[ Q_{j,1-\alpha_0}(\widehat{\boldsymbol{a}}_j^\beta) \pm z_{\alpha/2} \frac{\widehat{\sigma}(Q_{j,1-\alpha_0})}{N} \right]$$
and
$$IC_\alpha(Q_{j,1-\alpha_0}) = 	\left[Q_{j,1-\alpha_0}(\widehat{\boldsymbol{a}}_j^\beta) \exp\left(-z_{\alpha/2} \frac{\sigma(Q_{j,1-\alpha_0})}{Q_{1-\alpha_0}(\widehat{\boldsymbol{a}}_j^\beta)}\right), 
Q_{j,1-\alpha_0}(\widehat{\boldsymbol{a}}_j^\beta) \exp\left(z_{\alpha/2} \frac{\sigma(Q_{j,1-\alpha_0})}{Q_{j,1-\alpha_0}(\widehat{\boldsymbol{a}}_j^\beta)}\right)\right],$$
where $z_{\alpha}$ is as defined before.


\section{Influence function analysis}

The robustness of an estimator can be theoretically evaluated through its influence function (IF). Intuitively, the IF of an estimator measures how sensitive the estimator is to small changes in the observed data. 
It can be visualized as a measure of the pull that a single data point exerts on the estimator. Thus, robust estimators are expected to have a bounded IF, meaning that perturbations of a single data point have a bounded influence on the estimator.
The IF of an estimator is defined in terms of its statistical functional, defining the estimator for a general underlying distribution.
Let us consider $\widehat{\boldsymbol{a}}^\beta$ the MDPDE and $\boldsymbol{T}_\beta(G)$ its associated statistical functional. 
  That is, the MDPDE is defined as $\widehat{\boldsymbol{a}}^\beta = \boldsymbol{T}_\beta(G_n)$ 
 with $G_n$ the empirical distribution function of the lifetime estimated from the observed data. Furthermore, let us consider a contaminated lifetime distribution at a point perturbation $t_0,$ $G_\varepsilon = (1-\varepsilon)G + \varepsilon\Delta_{t_0},$ where $\varepsilon$ denotes the contamination proportion and $\Delta$ is a degenerate distribution at $t_0$. 
  Then, the IF of the MDPDE at a point perturbation $t_0$ and the theoretical lifetime distribution $F_{T}$ is given by
\begin{equation}\label{eq:IF}
	\text{IF}\left(t_0, \boldsymbol{T}_\beta, F_{T} \right) = \boldsymbol{J}_\beta^{-1}(\boldsymbol{a}_0) \boldsymbol{W}^T \boldsymbol{D}_{\boldsymbol{p}(\boldsymbol{a}_0)}^{\beta-1}\left(-\boldsymbol{p}(\boldsymbol{a}_0)+\delta_{t_0} \right),
\end{equation}
where $\boldsymbol{J}_\beta(\boldsymbol{a})$ is as defined in (\ref{eq:IF}), $\boldsymbol{W}$ and  $\boldsymbol{D}_{\boldsymbol{p}(\boldsymbol{a}_0)}$ are as defined in (\ref{eq:estimating}), $\boldsymbol{p}(\boldsymbol{a}_0)$ is the probability vector of the multinomial model with probabilities given in (\ref{eq:th}), and $\delta_{t_0}$ is an $LR+1$-dimensional degenerate probability vector with probability 1 at the cell containing the perturbation point $t_0.$ Indeed, the points of contamination on the multinomial model are given by multinomial samples 
with all components equal to $0$ except for a component, say $i_0,$ with $N$ failures.

For the step-stress ALT model with interval-censoring and competing risks, we should consider two forms of contamination: contamination in the covariates, including inspection times and stress levels, and outliers in the response variable given by the failure counts of each inspected cell. \cite{balakrishnan2022robust} studied the boundedness of the IF of the MDPDE for the  step-stress ALT model with interval-censoring and exponential lifetimes for contamination in the covariates, and established that the IF was bounded only for positives values of $\beta$. That is, the only non-robust estimator in the DPD family is the MLE, corresponding to $\beta=0$. Similar arguments can be used to discuss the boundedness of the IF of the MDPDE in the scenario of competing risks.
On the other hand, since the multinomial model has discrete support, the IF of the MDPDE at the outlier cells is bounded for any $\beta\geq 0,$ including the MLE. For discrete models where the IF is always bounded, the gross error sensitivity can be used for comparing the robustness of different estimators.
 We define the gross error sensitivity of the
functional $\boldsymbol{T}_\beta$ considering the contamination point at the $i_0-$cell as
\begin{equation}
	\gamma(i_0, \boldsymbol{T}_\beta, F_{T} )
	= \operatorname{sup}_{\{i_0=1,...,L+1\}} {\parallel IF\left(i_0, \boldsymbol{T}_\beta, F_{T} \right) \parallel}.
\end{equation}

\begin{figure}[htb]
 	\centering
 	\includegraphics[width = 10cm, height=5cm]{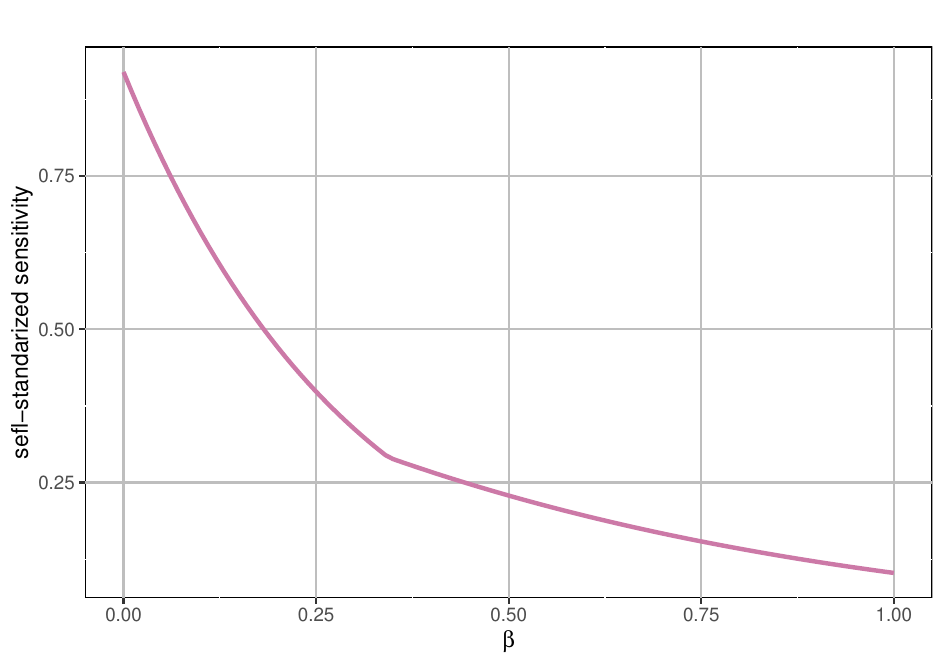}
 	\caption{Self-standardized sensitivity of the MDPDE with increasing DPD tuning parameter $\beta$}
 	\label{fig:sensitivity}
 \end{figure}

This maximum value, which describes the maximum bias on the
MDPDE over the neighborhood of the assumed model distribution, can be very large depending on the choice of the underlying distribution parameters, particularly for scale families.
It is natural to postulate that a measure of the robustness of parameter estimators should be invariant to scale transformations of individual parameter components at least-but the sensitivity defined above is not. To overcome this, \cite{hampel2011robust} measured the IF, which is an asymptotic bias, in
the metric given by the asymptotic covariance matrix of the estimator. That is, given the asymptotic variance-covariance matrix of the MDPDE defined in (\ref{eq:Sigma}), $\boldsymbol{\Sigma}(\boldsymbol{a}_0),$ the self-standardized sensitivity is defined
by

\begin{equation}
	\gamma^\ast(i_0, \boldsymbol{T}_\beta, F_{T} )
	= \operatorname{sup}_{\{i_0=1,...,L+1\}} {\parallel IF\left(i_0, \boldsymbol{T}_\beta, F_{T} \right)^T \boldsymbol{\Sigma}^{-1}(\boldsymbol{a}_0)IF\left(i_0, \boldsymbol{T}_\beta, F_{T} \right) \parallel}.
\end{equation}
 To illustrate the different sensitivities of the MDPDEs, Figure \ref{fig:sensitivity} presents, for increasing $\beta$s, the  self-standardized sensitivity of the estimators for a step-stress ALT experiment with two competing risks following exponential distributions with parameters $\theta_{1} = \exp(5-0.02x_1) $ and $\theta_{2} = \exp(6-0.04x_2).$ 
 The remaining components of the experimental design remain unchanged from the simulation scenario presented in Section \ref{sec:simulation}. It can be readily seen in Figure \ref{fig:sensitivity} that the standardized sensitivity decreases as $\beta$ increases, exhibiting the gain in robustness for larger values of $\beta.$
 


%

\section{Simulation study \label{sec:simulation}}

We evaluate the performance of the MDPDEs through  Monte Carlo simulation. 
For this purpose, we simulate a simple step-stress ALT with two stress levels, namely, $x_1= 35$ and $x_2=45,$ and two competing risks.
The test is performed for $75$ hours and the stress level is increased at $\tau_1= 45$ hours.
A total of $N=360$ devices are tested and a functional test is performed at times  $15,25,35,45,55,65$ and $75$ hours, where the number of failures attributed to each competing risk is recorded.

The lifetimes of the devices are generated from  exponential distributions with true parameters $\boldsymbol{a}_1 = (a_{01}, a_{11}) = (5,-0.02)$ for the first and $\boldsymbol{a}_2 = (a_{02}, a_{12}) = (6.2,-0.04)$ for the second competing risk.

To evaluate the robustness of the method, we introduce $\varepsilon\%$ of outlying observations coming from a distribution different from the underlying one. In particular, we consider a multinomial distribution with equal probabilities of failure within the third, fourth, fifth and sixth inspected intervals. That is, we increment the number of failures for the second and third inspected intervals.
We take the percentage of outliers as $\varepsilon = 0,0.05,0.1,0.2,0.3,0.4$ and $0.6.$

Figure \ref{fig:MSE} shows the RMSE on the estimation for different values of the tuning parameter $\beta.$ It can be seen that the performance of the MLE ($\beta=0$) rapidly worsens when introducing contamination in the sample, while robust estimators remain competitive even under high contamination proportions.

\subsection{Mean squared error of the model parameter estimates}

We first evaluate the performance of the MDPDE with different values of $\beta$ under increasing contamination percentages. Since the value of the true parameter comprises different magnitudes, we compute the  mean squared error in the estimation averaged in $R=1000$ simulation of each parameter $\widehat{a}_{ij}$ $i=1,2$ and $j=1,...,R,$ as follows:
$$\text{MSE} (\widehat{a}_{ij})= \frac{1}{R} \sum_{r=1}^R  (\widehat{a}_{ij}-a_{ij})^2.$$

\begin{figure}
	\centering
	\begin{subfigure}{0.4\textwidth}
		\includegraphics[scale=0.5]{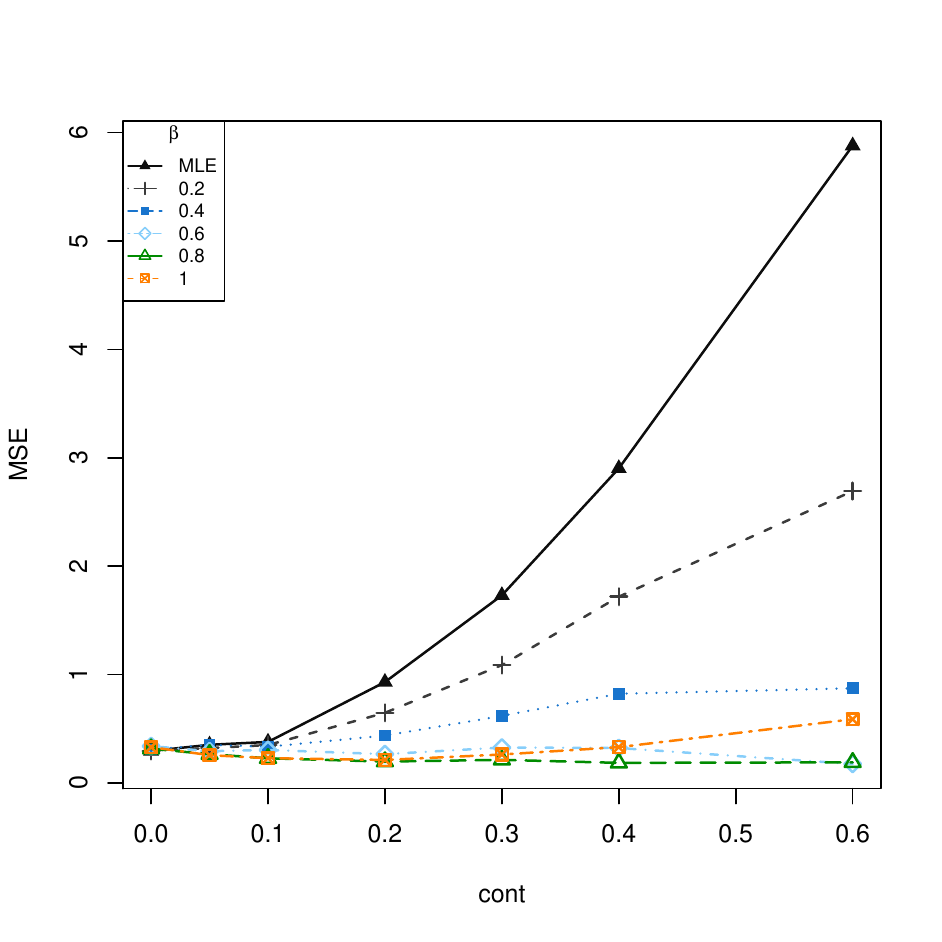}
		\caption{RMSE($\widehat{a}_{01}$)}
	\end{subfigure}
	\begin{subfigure}{0.4\textwidth}
		\includegraphics[scale=0.5]{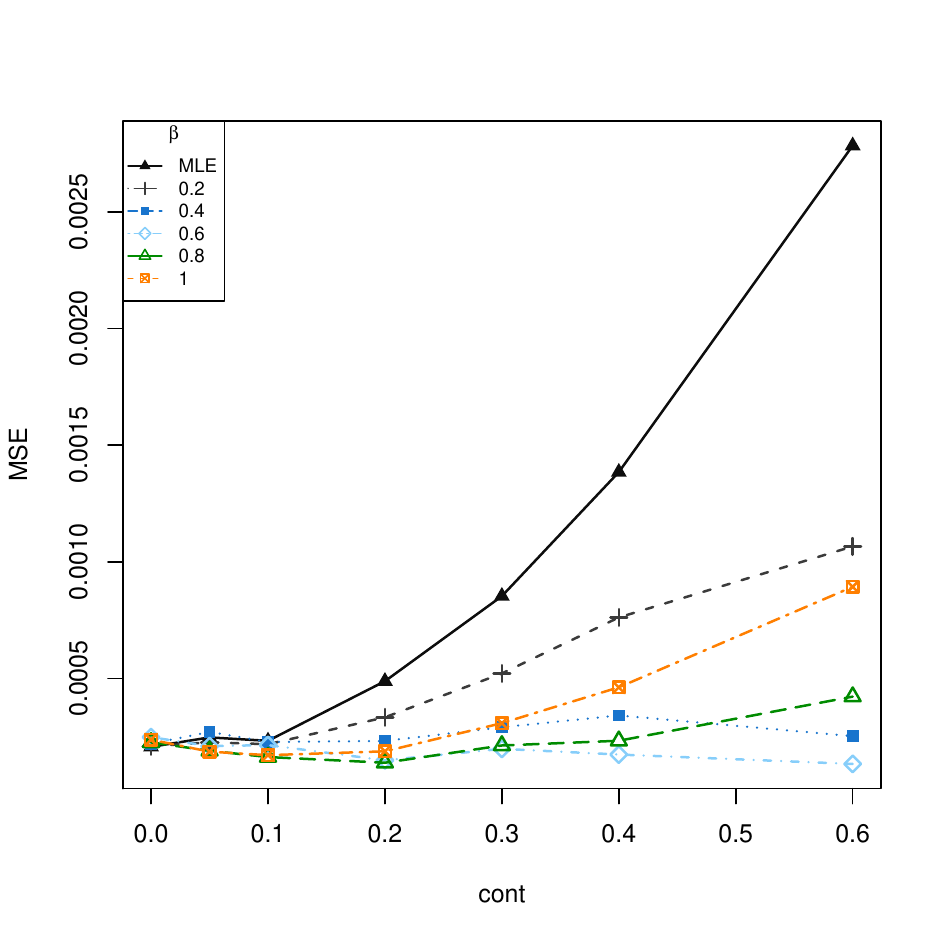}
		\caption{RMSE($\widehat{a}_{11}$)}
	\end{subfigure}	\begin{subfigure}{0.4\textwidth}
		\includegraphics[scale=0.5]{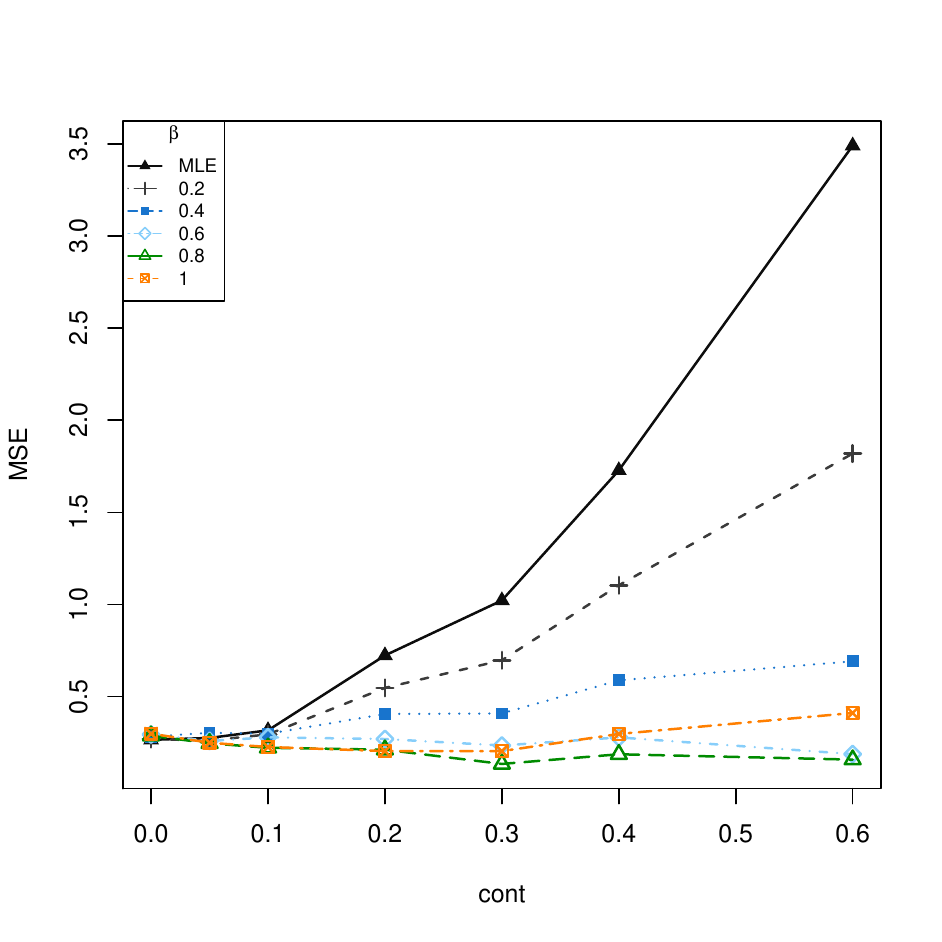}
		\caption{RMSE($\widehat{a}_{02}$)}
	\end{subfigure}
	\begin{subfigure}{0.4\textwidth}
		\includegraphics[scale=0.5]{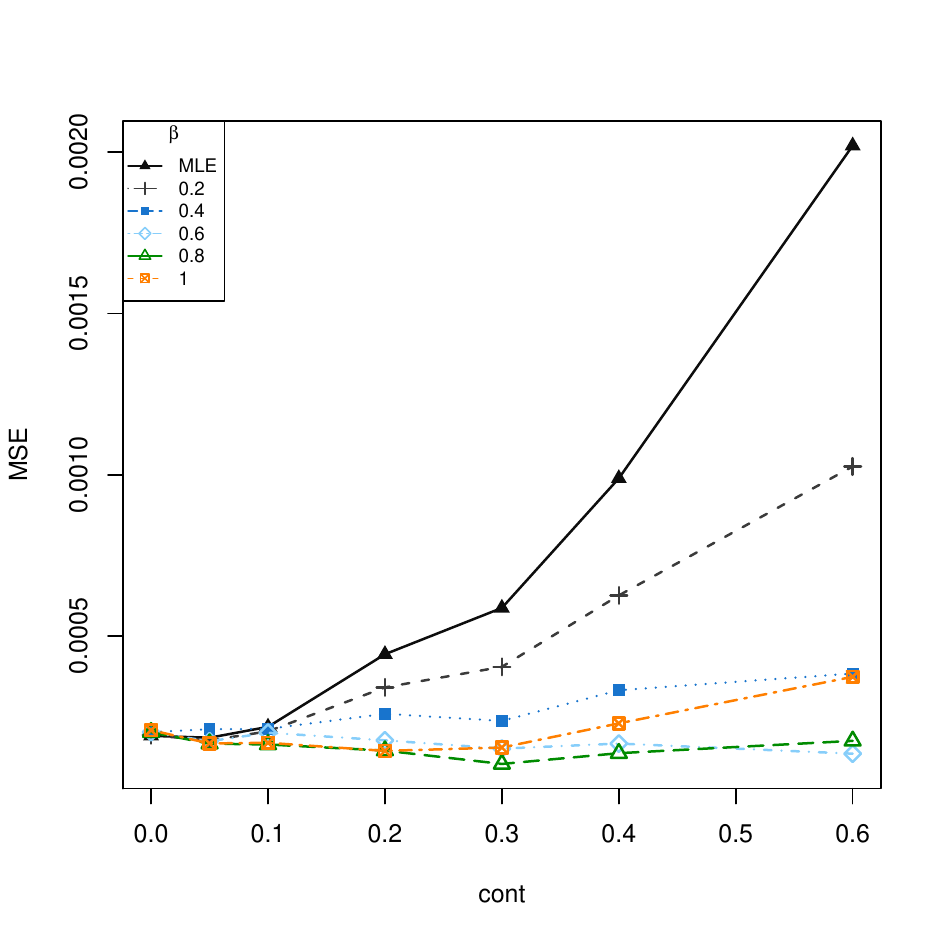}
		\caption{RMSE($\widehat{a}_{12}$)}
	\end{subfigure}
	\caption{Mean squared error of the MDPDEs with different values of the parameters $\boldsymbol{\theta}$}
	\label{fig:MSE}
\end{figure}

Figure \ref{fig:MSE} presents the averaged MSEs of the MDPDE under an increasing percentage of outlier observations for different values of the tuning parameter $\beta = 0, 0.2,0.4,0.6,0.8,1.$ From the plot, it is evident that the estimation error  committed by the MLE explodes from low contamination levels, unlike the more robust estimators with large $\beta$s. The use of positive values of $\beta,$ even small and moderate ones, has a significant impact on the robustness of the estimators.
Note that the improvement on robustness is shown for the estimation of all model parameters, and therefore it is transferred to the estimation of marginal distributions of the competing risks.
In contrast, all MDPDEs perform quite similarly in the absence of contamination,  with the MLE being slightly more efficient than the others. 
From the above, we suggest moderate values of $\beta$ (over 0.4) to achieve a suitable compromise between robustness and efficiency.

\subsection{Estimation of lifetime characteristics}

As discussed in Sections \ref{sec:char} and \ref{sec:specificchar}, in many reliability analysis, one is often interested in estimating some lifetime characteristics rather than the complete distribution function of the lifetime. We now evaluate the error made by the different MDPDEs in the estimation of the mean lifetime at normal operating conditions, reliability at a fixed mission time, and quantiles of the distribution. We also analyze their corresponding direct asymptotic, transformed and bootstrap $95\%$ confidence intervals in terms of coverage probability and  average width.

\begin{figure}[htb]
	\centering
	\begin{subfigure}{0.3\textwidth}
		\includegraphics[scale=0.35]{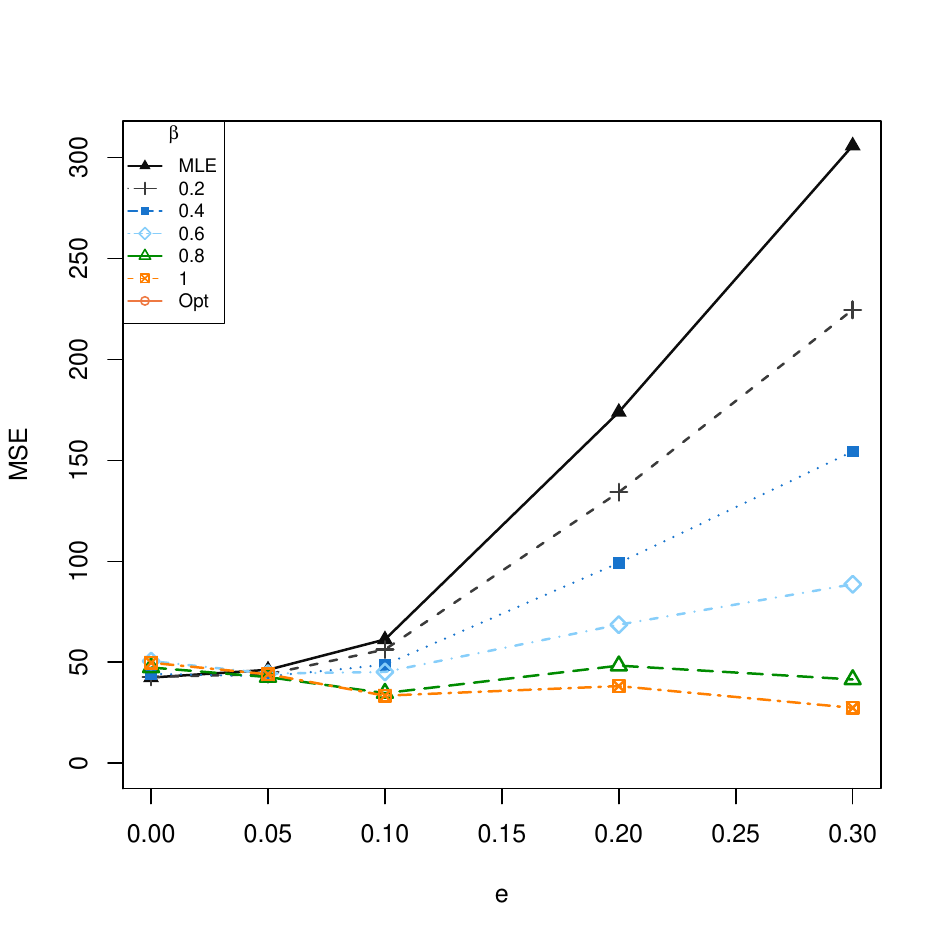}
		\caption{MTTF}
	\end{subfigure}
	\begin{subfigure}{0.3\textwidth}
		\includegraphics[scale=0.35]{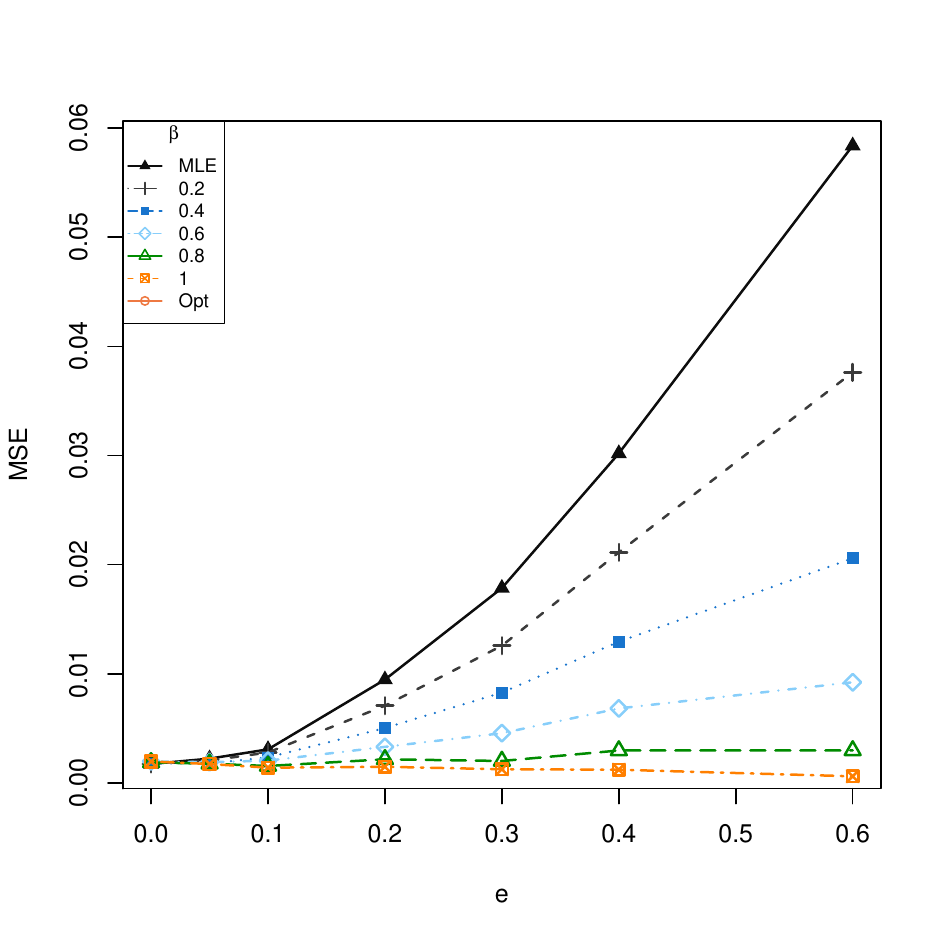}
		\caption{Reliability at $t_0=50$}
	\end{subfigure}	
	\begin{subfigure}{0.3\textwidth}
		\includegraphics[scale=0.35]{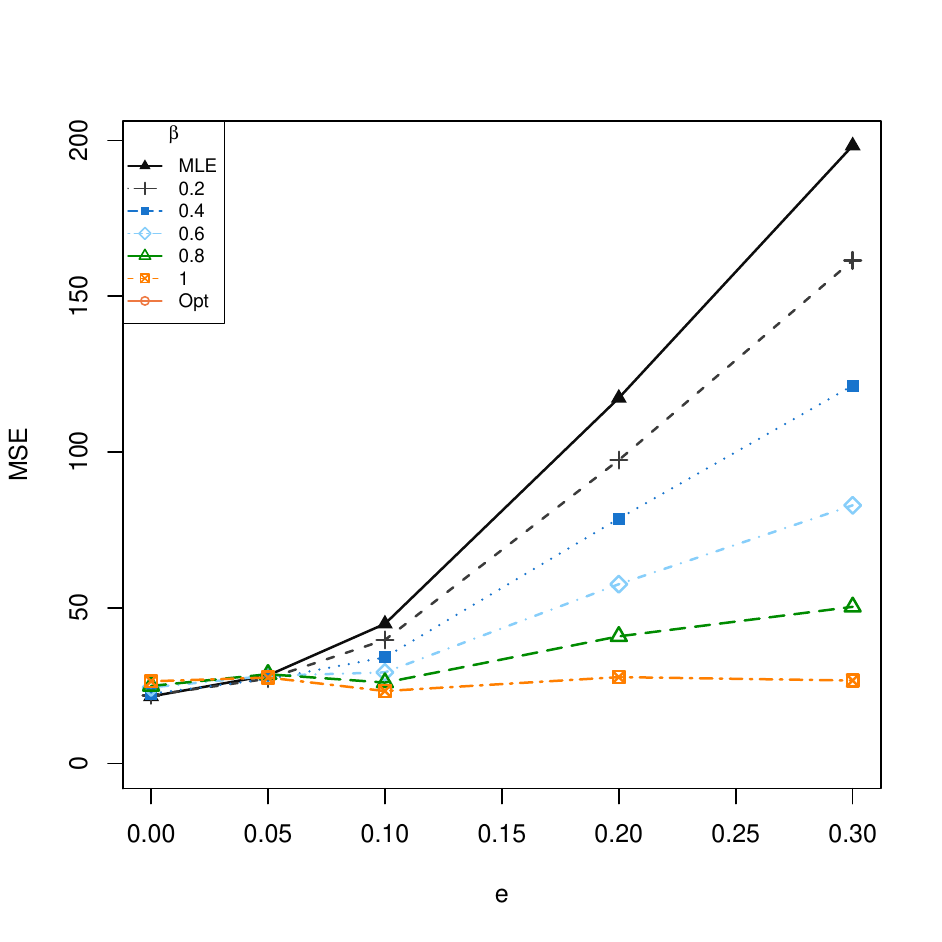}
		\caption{Median}
	\end{subfigure}
	\caption{Mean squared error (MSE) of the MDPDE with different values of $\beta$ for estimating the mean lifetime to failure, the relibility of the devices al $t_0=50,$ and the median of the distribution. }
	\label{fig:MSEchar}
\end{figure}

Figure \ref{fig:MSEchar} shows the performance of the MDPDEs with different values of the tuning parameter $\beta$ when estimating the three lifetime characteristics of interest, namely, the mean, median and reliability at a mission time of $t_0 = 50$. The robustness of the estimators is clearly inherited by the corresponding lifetime estimates, and again all estimates based on the MLE become worse rapidly when  outlying observations are introduced.

On the other hand, Table \ref{table:cis} shows the performance of the confidence intervals for the mean lifetime to failure in terms of coverage and width for increasing contamination, obtained with direct asymptotic, transformed and bootstrap intervals. Results for the median and reliability are quite similar and we therefore omitted them for brevity.
 Several conclusions can be drawn from the presented results. Firstly, the results illustrate the advantage of using transformed intervals instead of direct asymptotic confidence intervals under any contamination scenario and tuning parameter. Secondly, the gain in robustness is evidenced in terms of coverage without a great increase in the width for asymptotic and bootstrap confidence intervals.
Also,  the bootstrap confidence intervals naturally worsen considerably when the initial estimate is far from the true model. 
Moreover, although the MDPDEs with positive values of the tuning parameter present wider confidence intervals than the MLE, the loss in efficiency more than compensated by their robustness.


\begin{table}[ht]
	\centering
	\begin{tabular}{r|rr|rr|rr}
		\hline
		\multicolumn{7}{c}{$\varepsilon=0$}\\
		\hline
		& \multicolumn{2}{c}{Direct} & \multicolumn{2}{c}{Transformed} & \multicolumn{2}{c}{BCa}\\
		& Cov. & Width & Cov. & Width &  Cov. & Width  \\ 
		\hline
		$\beta=0$ & 100$\%$ & 39.30 & 100$\%$ & 40.05 & 86$\%$ & 22.50 \\ 
		$\beta=0.2$ & 100$\%$ & 47.41 & 100$\%$ & 48.73 & 87$\%$ & 21.75 \\ 
		$\beta=0.4$ & 100$\%$ & 57.30 & 100$\%$ & 59.64 & 90$\%$ & 23.00 \\ 
		$\beta=0.6$ & 100$\%$ & 70.14 & 100$\%$ & 74.44 & 84$\%$ & 22.84 \\ 
		$\beta=0.8$& 100$\%$ & 85.16 & 100$\%$ & 93.04 & 86$\%$ & 23.33 \\ 
		$\beta=1$ & 100$\%$ & 104.47 & 100$\%$ & 119.34 & 88$\%$ & 23.67 \\ 
		\hline
	\multicolumn{7}{c}{$\varepsilon=5\%$}\\
		\hline
		$\beta=0$ & 96$\%$ & 37.04 & 100$\%$ & 37.75 & 72$\%$ & 20.95 \\ 
		$\beta=0.2$ & 100$\%$ & 45.12 & 100$\%$ & 46.39 & 79$\%$ & 20.73 \\ 
		$\beta=0.4$ & 100$\%$ & 55.95 & 100$\%$ & 58.28 & 84$\%$ & 21.87 \\ 
		$\beta=0.6$ & 100$\%$ & 68.15 & 100$\%$ & 72.40 & 78$\%$ & 22.18 \\ 
		$\beta=0.8$ & 100$\%$ & 84.54 & 100$\%$ & 92.58 & 85$\%$ & 23.12 \\ 
		$\beta=1$ & 100$\%$ & 105.57 & 100$\%$ & 121.21 & 86$\%$ & 23.51 \\ 
		\hline
	\multicolumn{7}{c}{$\varepsilon=10\%$}\\
		\hline
		$\beta=0$ & 91$\%$ & 35.20 & 96$\%$ & 35.89 & 61$\%$ & 19.77 \\ 
		$\beta=0.2$ & 97$\%$ & 43.17 & 100$\%$ & 44.40 & 69$\%$ & 19.58 \\ 
		$\beta=0.4$ & 100$\%$ & 53.87 & 100$\%$ & 56.16 & 75$\%$ & 20.59 \\ 
		$\beta=0.6$ & 100$\%$ & 67.36 & 100$\%$ & 71.68 & 81$\%$ & 21.15 \\ 
		$\beta=0.8$ & 100$\%$ & 84.37 & 100$\%$ & 92.66 & 88$\%$ & 22.64 \\ 
		$\beta=1$ & 100$\%$ & 106.93 & 100$\%$ & 123.46 & 92$\%$ & 23.29 \\ 
		\hline
	\multicolumn{7}{c}{$\varepsilon=20\%$}\\
		\hline
		$\beta=0$ & 69$\%$ & 31.57 & 75$\%$ & 32.21 & 31$\%$ & 17.53 \\ 
		$\beta=0.2$ & 85$\%$ & 39.48 & 96$\%$ & 40.65 & 39$\%$ & 17.62 \\ 
		$\beta=0.4$ & 99$\%$ & 50.51 & 100$\%$ & 52.74 & 49$\%$ & 18.58 \\ 
		$\beta=0.6$ & 100$\%$ & 64.25 & 100$\%$ & 68.54 & 64$\%$ & 19.69 \\ 
		$\beta=0.8$ & 100$\%$ & 82.88 & 100$\%$ & 91.48 & 76$\%$ & 21.20 \\ 
		$\beta=1$ & 100$\%$ & 108.73 & 100$\%$ & 126.78 & 87$\%$ & 22.63 \\ 
		\hline
	\multicolumn{7}{c}{$\varepsilon=30\%$}\\
		\hline
		& \multicolumn{2}{c}{Direct} & \multicolumn{2}{c}{Transformed} & \multicolumn{2}{c}{BCa}\\
		\hline
		$\beta=0$ & 34$\%$ & 29.32 & 59$\%$ & 29.94 & 3$\%$ & 15.95 \\ 
		$\beta=0.2$ & 75$\%$ & 37.55 & 90$\%$ & 38.71 & 11$\%$ & 16.04 \\ 
		$\beta=0.4$ & 97$\%$ & 48.75 & 100$\%$ & 51$\%$ & 26$\%$ & 17.56 \\ 
		$\beta=0.6$ & 100$\%$ & 64.41 & 100$\%$ & 68.99 & 58$\%$ & 19.40 \\ 
		$\beta=0.8$ & 100$\%$ & 85.48 & 100$\%$ & 95.06 & 79$\%$ & 20.66 \\ 
		$\beta=1$ & 100$\%$ & 114.65 & 100$\%$ & 135.50 & 96$\%$ & 22.50 \\ 
		\hline
	\end{tabular}
\caption{Coverage (Cov) and width of the estimated confidence intervals for the mean lifetime to failure using direct asymptotic  and transformed confidence intervals as well as bootstrap bias corrected (BCa) intervals for different values of the tuning parameter $\beta$ and different proportions of contaminated observations $\varepsilon.$}
\label{table:cis} 
\end{table}



\section{Analysis of Electronic Device Data}

\cite{han2014inference} conducted a simple step-stress test under time constraints to evaluate the reliability characteristics of a solar lighting device.  The
dataset consists of a total of 31 failure times  from the initial sample size of $N= 35$ prototypes.
The device was found to be susceptible to two main failure modes: capacitor failure and controller failure. The test focused on the influence of temperature, which served as the stress factor. The temperature was varied within a range of 293K to 353K, with the normal operating temperature set at 293K.
The stress level was changed at a time point of $\tau_1 = 500$ hours, while the end of the test was set as $\tau_2$ = 600  hundred hours.
Exact failure times were analyzed in their study. However, for illustrative  purposes, we will interval-censor the data at inspection times $IT = 200,400,500,525,550$ and $600$ hours.
The interval-censored data so obtained is presented in Table \ref{table:realdata}.
\begin{table}
	\centering
	\begin{tabular}{ccc}
		IT &  Number of capacitor failures &  Number of  controller failures\\
		\hline
		200 & 2 & 5 \\
		300 & 6 & 0 \\
		500 &  2 & 1 \\
		525 & 2 & 6 \\
		550 & 2 & 4\\
		600 & 1 & 0\\
		\hline
	\end{tabular}
\caption{Interval-monitored solar lighting device data }
\label{table:realdata}
\end{table}

We assume that the lifetime distributions of the device's failure due to each risk factor, namely, capacitor failure and controller failure, are independent and follow an exponential distribution at any constant temperature.
The Arrhenius law relates the exponential scale parameter $\theta_i$ to the experienced temperature $T_i$ as
$$\theta_i = \exp\left(-\frac{E_a}{K}\left(\frac{1}{T_0}- \frac{1}{T_i}\right) \right) \theta_0,$$
where $\theta_0$ is the scale parameter of the lifetime distribution under normal operating temperature $T_0,$
$K =  8.36 \times 10^{-5}e V/^\circ K$ is Boltzmann’s constant, and $E_a$ is the activation energy.
Then, re-parametrizing the stress level as $$x_i=-\frac{1}{K}\left[\frac{1}{T_0}-\frac{1}{T_i}\right],$$
the Arrhenius  equation can be rewritten as a log-linear relation of the form
 $$\theta_i = \exp\left(\log(\theta_0) +E_a x_i \right).$$  
 %
In our example, we apply the above transformation to the temperatures $T_0=293$K and $T_1=353$K, obtaining the stress levels $x_0 = 0$ and $x_1=1.$
Table \ref{table:realdata-characteristics} presents the MDPDE of the mean lifetime to failure, median and reliability at mission time $t_0=4$ of the electronic devices' lifetime with different values of $\beta,$ jointly with their corresponding direct asymptotic, transformed and bootstrap confidence intervals. MDPDEs with large values of $\beta$ estimate larger mean and median lifetimes and consequently find the devices to be more reliable.

For the electronic devices dataset, direct confidence intervals for the three lifetime characteristics of interest had to be truncated and are therefore not very informative, especially for the reliability of the devices.
On the other hand, employing transformed confidence intervals solve the truncation issue and provide more informative intervals, yet still quite wide. Notably, for large values of the tuning parameter $\beta,$ the estimated  standard error of the MDPDE is quite large, resulting in too wide intervals.
On the other hand, bootstrap techniques offer an appealing alternative for small sample sizes, such as the one under study, with narrower confidence intervals. The wider confidence intervals and a considerably high coverage probability are due to the fact that the data involves a somewhat small sample size with a number of inspection points thus making the number of failures in each time interval smaller. This problem will be alleviated when the number of test units is large.

From the results, we observe that
moderate and small values of $\beta$ may be better to used. Moreover, transformed confidence intervals based on the asymptotic distribution of the estimators provides more insightful intervals, although bootstrap techniques provide narrower intervals which may be certainly preferable in case of small sample sizes.

\begin{table}[ht]
	\centering
	\begin{tabular}{ccccc}
		\hline
		\multicolumn{5}{c}{Mean lifetime to failure}\\
		\hline
		& Estimate & Direct CI & Transformed CI & Bootstrap CI \\ 
		\hline
		$\beta=0$ & 7.81 & [0.00, 21.72] & [1.32, 46.34]  & [5.32, 11.19]\\ 
		$\beta=0.2$ & 8.06 & [0.00, 26.64] & [0.81, 80.73] & [5.54, 11.73] \\ 
		$\beta=0.4$  & 8.30 & [0.00, 33.06] & [0.42, 163.81] & [5.44, 11.47]\\ 
		$\beta=0.6$ & 8.51 & [0.00, 41.34] & [0.18, 402.64] & [5.50, 12.57]\\ 
		$\beta=0.8$ & 8.69 & [0.00, 51.91] & [0.06, 1258.86] & [5.50, 12.57]\\ 
		$\beta=1$ & 8.82 & [0.00, 65.22] & [0.01, 5296.84] & [5.61, 14.14]\\ 
%
		\hline
		\multicolumn{5}{c}{Reliability}\\
		\hline
		& Estimate & Direct CI & Transformed CI & Bootstrap CI \\  
		\hline
		$\beta=0$ & 0.60 & [0.05, 1.00]& [0.13, 0.94]&[0.43, 0.72] \\ 
		$\beta=0.2$ & 0.61 & [0.00, 1.00] & [0.08, 0.97] & [0.63, 0.71] \\ 
		$\beta=0.4$ & 0.62 & [0.00, 1.00] & [0.04, 0.99]& [0.50, 0.74] \\ 
		$\beta=0.6$ & 0.63 & [0.00, 1.00] & [0.01, 1.00]& [0.50, 0.74]\\ 
		$\beta=0.8$ & 0.63 & [0.00, 1.00] & [0.00, 1.00] & [0.49, 0.73]\\ 
		$\beta=1$ & 0.64 & [0.00, 1.00] & [0.00, 1.00] & [0.52, 0.75] \\ 
		\hline
%
	\multicolumn{5}{c}{Median}\\
	\hline
	& Estimate & Direct CI & Transformed CI & Bootstrap CI \\ 
	\hline
	$\beta=0$ & 5.41 & [0.00, 15.05] & [0.91, 32.12]&[3.67,  8.15]\\ 
	$\beta=0.2$& 5.59 & [0.00, 18.47] & [0.56, 55.96] &  [3.35,  7.52] \\ 
	$\beta=0.4$& 5.76 & [0.00, 22.92] & [0.29, 113.55] & [3.92,  8.15]\\ 
	$\beta=0.6$ &  5.90  & [0.00, 41.34] & [0.18, 402.64] & [3.90,  9.19]\\ 
	$\beta=0.8$ & 6.02 & [0.00, 35.98] & [0.04, 872.58] &[4.14,  8.73] \\ 
	$\beta=1$  & 6.11 & [0.00, 45.21] & [0.01, 3671.49]  &[3.71, 11.34]\\ 
		\hline
	\end{tabular}
\caption{Estimates and confidence intervals of the mean, median, and reliability at mission time $t_0$ under normal operating conditions $T_0=293$ based on the MDPDEs with different values of the tuning parameter $\beta.$ }
\label{table:realdata-characteristics}
\end{table}

\section{Conclusions}

Step-stress accelerated life-tests enable the estimation of lifetime characteristics for highly reliable products. 
In some applications, the lifetime status of the units under test can not be continuously monitored  and so interval-censored data may arise.
Additionally, for such devices, there may exist several causes of failure, giving rise to competing risks.
In this paper, both robust and non-robust inferential methods for an interval-monitoring step-stress model under competing risks, are presented. Theoretical asymptotic properties of the MDPDEs, including the classical MLE, have been studied and approximate confidence intervals for the model's parameters have been presented.
Furthermore, point estimation and confidence intervals for some lifetime characteristics and cause-specific lifetime characteristics, namely, mean lifetime to failure, reliability at a mission time and distribution quantiles, have been developed. Transformed and bootstrap confidence intervals have been proposed as two alternatives to the direct asymptotic intervals. Their performance have been empirically compared, showing the pros and cons of each approach.
Finally, the robustness of the estimators has been analyzed theoretically and empirically through an extensive simulation study. Finally,
the usefulness of the inferential methods with a real-life testing is illustrated on a real dataset regarding the reliability of electronic devices.
It will naturally be of interest to consider more general lifetime distributions for the competing causes such as Weibull and Gamma, and then develop analogous results. We are working on this problem.

\bibliographystyle{abbrvnat}
\bibliography{bibliography}

\section*{Acknowledgements}
This work was supported by the Spanish Grant PID2021-124933NB-I00  and the Natural Sciences and Engineering Research Council of Canada (of the first author) through an Individual Discovery Grant (No. 20013416).
M. Jaenada and L. Pardo are members of the Interdisciplinary Mathematics Institute (IMI).

\appendix
\section{Computation of derivatives of the probability of failure \label{AppA}}

We derive here explicit expressions for the derivatives with respect to the parameters $a_{0k}$ and $a_{1k}$ of the probability of failure in the $l$-th interval due to risk $j$, for any $j, k = 1,...,R$ and $l = 1,...,L.$
Recall that the risk imposed on a test unit due to risk factor $j$ is given by
$$\pi_{ij}(\boldsymbol{a}) = \frac{\exp(-a_{0j}-a_{1j}x_i)}{\sum_{j=1}^R \exp(-a_{0j}-a_{1j}x_i)}.$$
Let us define an auxiliary function 
$$S(\boldsymbol{a}, IT_l) = \exp\left(-\sum_{j=1}^{R} \frac{IT_{l}+h_j^{(i)}}{\exp(a_{0j}+a_{1j}x_i)}\right)$$
with
\begin{equation*} 
	h^{(1)}_j = 0 \hspace{0.3cm} \text{and} \hspace{0.3cm} h^{(2)}_j = \tau_1\left(\frac{\exp(a_{0j}+a_{1j}x_2)}{\exp(a_{0j}+a_{1j}x_1)}- 1\right)
\end{equation*}
as defined in (\ref{eq:shiftingtime2}). 
Then, the probability  of failure at the $l$-th interval due to risk $j$ can be rewritten as
\begin{equation}\label{eq:simplyptheta}
		p_{lj}(\boldsymbol{a})  =  
		\pi_{ij}(\boldsymbol{a})
		 \left( S(\boldsymbol{a}, IT_{l-1}) - S(\boldsymbol{a}, IT_l)  \right), \hspace{0.3cm} \tau_{i-1} \leq IT_l < \tau_i.
\end{equation}
From Equation (\ref{eq:simplyptheta}), it is natural to define the following quantities:
\begin{equation}\label{eq:simplyptheta2}
	R_{lj}^k(\boldsymbol{a})  =  
	\pi_{ij}(\boldsymbol{a}) S(\boldsymbol{a}, IT_{l}) \hspace{0.3cm} \text{with} \tau_{i-1} \leq IT_l < \tau_i
\end{equation}
and then taking derivatives in (\ref{eq:simplyptheta}), we get
\begin{equation}
	\begin{aligned}
		\frac{\partial p_{lj}(\boldsymbol{a})}{\partial a_{0k}} &=  \frac{\partial \pi_{ij}(\boldsymbol{a})}{\partial a_{0k}}\left( S(\boldsymbol{a}, IT_{l-1}) - S(\boldsymbol{a}, IT_l)  \right) +\pi_{ij}(\boldsymbol{a})  \left(\frac{\partial S(\boldsymbol{a}, IT_{l-1})}{\partial a_{0k}} - \frac{\partial S(\boldsymbol{a}, IT_l)}{\partial a_{0k}}  \right), \\
		& = \frac{\partial R_{l-1,j}(\boldsymbol{a})}{\partial a_{0k}} - \frac{\partial R_{lj}(\boldsymbol{a})}{\partial a_{0k}}\\
		\frac{\partial p_{lj}(\boldsymbol{a})}{\partial a_{1k}} &=  \frac{\partial \pi_{ij}(\boldsymbol{a})}{\partial a_{1k}}\left( S(\boldsymbol{a}, IT_{l-1}) - S(\boldsymbol{a}, IT_l)  \right) +\pi_{ij}(\boldsymbol{a})  \left(\frac{\partial S(\boldsymbol{a}, IT_{l-1})}{\partial a_{1k}} - \frac{\partial S(\boldsymbol{a}, IT_l)}{\partial a_{1k}}  \right)\\
		& = \frac{\partial R_{l-1,j}(\boldsymbol{a})}{\partial a_{1k}} - \frac{\partial R_{lj}(\boldsymbol{a})}{\partial a_{1k}}.
	\end{aligned}
\end{equation}

We will compute the derivatives for the quantities $w_{lj}.$
First, note that the shifting times $h^{(i)}_j$ depends on the model parameters and its derivatives are given by
\begin{equation}
	\begin{aligned}
		\frac{\partial h^{(i)}_j}{\partial a_{0k}} 
		& = \frac{\partial h^{(1)}_j}{\partial a_{1k}}  = 0,\\
		\frac{\partial h^{(2)}_j}{\partial a_{1k}} & =
		\begin{cases}
			0 & k \neq j \\
			\tau_1 \frac{\exp(a_{0j}+a_{1j}x_2)}{\exp(a_{0j}+a_{1j}x_1)} (x_2-x_1), & k=j
		\end{cases}
	\end{aligned}
\end{equation}
and that the derivatives of the relative risks satisfy the equations
\begin{equation*}
	\frac{\partial \pi_{ij}(\boldsymbol{a})}{\partial a_{0k}} = \begin{cases}
	  \pi_{ij}(\boldsymbol{a}) \pi_{ik}(\boldsymbol{a}) & k \neq j,\\
	  \pi_{ij}(\boldsymbol{a})(\pi_{ij}(\boldsymbol{a})- 1), & k = j.
	\end{cases}
\hspace{0.5cm}
	\frac{\partial \pi_{ij}(\boldsymbol{a})}{\partial a_{1k}} = \begin{cases}
		 \pi_{ij}(\boldsymbol{a}) \pi_{ik}(\boldsymbol{a})x_i, & k \neq j,\\
		\pi_{ij}(\boldsymbol{a})(\pi_{ij}(\boldsymbol{a})- 1)x_i, & k = j.
	\end{cases}
\end{equation*}
Now, let us denote
\begin{equation}\label{eq:shiftingtime_ast}
	h_j^{\ast (i)} = \tau_1\frac{\exp(a_{0j}+a_{1j}x_2)}{\exp(a_{0j}+a_{1j}x_1)} (x_2-x_1).
\end{equation}
Then, taking derivatives in the auxiliary function $S(\boldsymbol{a}, IT_{l}),$ we get
\begin{equation*}
	\begin{aligned}
		\frac{\partial S(\boldsymbol{a}, IT_{l})}{\partial a_{0k}} &= 
		S(\boldsymbol{a}, IT_{l})\left(\frac{IT_l + h_k^{(i)}}{\exp(a_{0k}+a_{1k}x_i)}\right),\\
		\frac{\partial S(\boldsymbol{a}, IT_{l})}{\partial a_{1k}} &= 
		S(\boldsymbol{a}, IT_{l})\left(\frac{-h_k^{\ast (i)}+(IT_l + h_k^{(i)})x_i}{\exp(a_{0k}+a_{1k}x_i)}\right).
	\end{aligned}
\end{equation*}

So, we can compute the derivatives of $R_{lj}(\boldsymbol{a})$ in Equation (\ref{eq:simplyptheta2}) as 
\begin{equation*}
	\begin{aligned}
		\frac{\partial R_{lj}(\boldsymbol{a})}{\partial a_{0k}}= &= 
		\begin{cases}
			\pi_{ij}S(\boldsymbol{a}, IT_{l})\left(\pi_{ik}+\frac{IT_l + h_k^{(i)}}{\exp(a_{0k}+a_{1k}x_i)}\right), & j \neq k, \\
			\pi_{ij}S(\boldsymbol{a}, IT_{l})\left(\pi_{ij}-1+\frac{IT_l + h_k^{(i)}}{\exp(a_{0k}+a_{1k}x_i)}\right),  & j \neq k,\\
		\end{cases}\\
		\frac{\partial R_{lj}(\boldsymbol{a})}{\partial a_{1k}} &= 
		\begin{cases}
			\pi_{ij}S(\boldsymbol{a}, IT_{l})\left(\pi_{ik}x_i+\frac{-h_k^{\ast(i)} + (IT_l + h_k^{(i)})x_i}{\exp(a_{0k}+a_{1k}x_i)}\right), & j \neq k, \\
			\pi_{ij}S(\boldsymbol{a}, IT_{l})\left(\pi_{ij}x_i-x_i+\frac{-h_k^{\ast(i)} + (IT_l + h_k^{(i)})x_i}{\exp(a_{0k}+a_{1k}x_i)}\right), & j \neq k,\\
		\end{cases}
	\end{aligned}
\end{equation*}
and finally we can write explicit expressions for the derivatives of the probability of failure as
\begin{equation}
	\begin{aligned}
		\frac{\partial p_{lj}(\boldsymbol{a})}{\partial a_{0k}}&= \frac{\partial R_{l-1j}(\boldsymbol{a})}{\partial a_{0k}} - \frac{\partial R_{lj}(\boldsymbol{a})}{\partial a_{0k}},\\
		\frac{\partial p_{lj}(\boldsymbol{a})}{\partial a_{1k}}&= \frac{\partial R_{l-1j}(\boldsymbol{a})}{\partial a_{1k}} - \frac{\partial R_{lj}(\boldsymbol{a})}{\partial a_{1k}}.
	\end{aligned}
\end{equation}

\section{Computation of lifetime characteristics of devices}

We first derive the mean of the lifetime to failure $T,$ under normal operating conditions, with PDF given by
$$
	f_T(t) 
	=
	\sum_{j=1}^R \frac{1}{\theta_{0j}} \exp\left(- \sum_{j=1}^R \frac{t}{\theta_{0j}}\right)  \hspace{0.3cm} t>0
 $$
Then,
\begin{equation*}
	\begin{aligned}
		E[T] &= \int_{0}^\infty t \left(\sum_{j=1}^R \frac{1}{\theta_{0j}} \exp\left(- \sum_{j=1}^R \frac{t}{\theta_{0j}}\right) \right)dt =\sum_{j=1}^R \frac{1}{\theta_{0j}} \int_{0}^\infty   t \exp\left(- \sum_{j=1}^R \frac{t}{\theta_{0j}}\right) dt\\
		&=\sum_{j=1}^R \frac{\theta_{0j}^{-1}}{\sum_{j=1}^R \theta_{0j}^{-1}}  \int_{0}^\infty   \exp\left(- \sum_{j=1}^R \frac{t}{\theta_{0j}}\right)dt
		= \frac{1}{\sum_{j=1}^R \theta_{0j}^{-1}} .
	\end{aligned}
\end{equation*}

\end{document}